\documentclass[12pt]{article}
\usepackage{epsfig}
\usepackage{amsfonts}
\usepackage{amssymb}
\usepackage{amsmath}
\usepackage{amsthm}
\usepackage{latexsym}
\usepackage{graphicx}
\usepackage{bm}
\usepackage{tabularx}
\usepackage{booktabs} 
\usepackage{enumerate}
\usepackage{caption}
\usepackage[titletoc]{appendix}
\usepackage[numbers]{natbib}
\usepackage{bbm}
\usepackage{url}
\usepackage{subfig}
\usepackage{mathtools}

\catcode`\@=11

\newcommand{\shortdot}[1]{\raisebox{-0.4pt}{$\stackrel{\bullet}{#1}$}}
\newcommand{\updot}[1]{\raisebox{0.9pt}{$\stackrel{\bullet}{#1}$}} %%
\newcommand{\dshortdot}[1]{\raisebox{-0.4pt}{$\stackrel{\bullet \bullet}{#1}$}}

\usepackage[pdftitle={Delayed Information},
  pdfauthor={Pender, Rand, Wesson},
  pdffitwindow=true,
  breaklinks=true,
  colorlinks=true,
  urlcolor=blue,
  linkcolor=red,
  citecolor=red,
  anchorcolor=red]{hyperref}

\theoremstyle{plain}
\newtheorem{theorem}{Theorem}[section]

\newtheorem{corollary}[theorem]{Corollary}

\theoremstyle{definition}

\theoremstyle{remark}
\newtheorem*{remark}{Remark}

\begin{document}
\title{An Asymptotic Analysis of Queues with Delayed Information and Time Varying Arrival Rates}
%\title{An Asymptotic Analysis of Queues with Delayed Information and Time Varying Arrival Rates}
\author{ 
  Jamol Pender \\ School of Operations Research and Information Engineering \\ Cornell University
\\ 228 Rhodes Hall, Ithaca, NY 14853 \\  jjp274@cornell.edu  \\ 
 \and  
Richard H. Rand \\ Sibley School of Mechanical and Aerospace Engineering \\ Department of Mathematics \\ Cornell University
\\ 535 Malott Hall, Ithaca, NY 14853 \\  rand@math.cornell.edu  \\ 
 \and  
Elizabeth Wesson \\ Department of Mathematics \\ Cornell University
\\ 582 Malott Hall, Ithaca, NY 14853 \\  enw27@cornell.edu  \\ 
 }

\maketitle
\begin{abstract}
Understanding how delayed information impacts queueing systems is an important area of research.  However, much of the current literature neglects one important feature of many queueing systems, namely non-stationary arrivals.  Non-stationary arrivals model the fact that customers tend to access services during certain times of the day and not at a constant rate.  In this paper, we analyze two two-dimensional deterministic fluid models that incorporate customer choice behavior based on delayed queue length information with time varying arrivals.  In the first model, customers receive queue length information that is delayed by a constant $\Delta$.  In the second model, customers receive information about the queue length through a moving average of the queue length where the moving average window is $\Delta$.  We analyze the impact of the time varying arrival rate and show using asymptotic analysis that the time varying arrival rate does not impact the critical delay unless the frequency of the time varying arrival rate is twice that of the critical delay.  When the frequency of the arrival rate is twice that of the critical delay, then the stability is enlarged by a wedge that is determined by the model parameters.  As a result, this problem allows us to combine the theory of nonlinear dynamics, parametric excitation, delays, and time varying queues together to provide insight on the impact of information in queueing systems.    

\end{abstract}

%**************************************************************************
%**************************************************************************

\section{Introduction} \label{sec_intro}

Understanding the impact of providing delayed information to customers in queueing systems is a very important problem in the queueing and engineering literature.  Many companies where customers are forced to wait in line often choose to provide their customers with waiting time or queue length information.  Consequently, the information that is provided can affect a customer's choice of using the service and joining the queue.  One common example of this communication between the service and customer is delay announcements.  Delay announcements commonly inform customers about the average waiting time to start service.  These announcements are not only important because they give the customer information about their wait, but also the announcements have the possibility of influencing the possibility that a customer will return to use the service again or remain in.  As a consequence, understanding the impact of providing queue length information to customers on customer choices and system operations, as well as the development of methods to support such announcements, has attracted the attention of the queueing systems community recently. 

Much of the research on providing queue length or waiting time information to customers focuses on the impact of delay announcements in call centers.  There is a vast literature on this subject, but our focus is quite different.  Work by \citet{ibrahim2008real, ibrahim2009real, ibrahim2011real, ibrahim2011wait} develops new real-time estimators for estimating delays in various queueing systems.  The work of \citet{armony2004customer, guo2007analysis, hassin2007information, armony2009impact, guo2009impacts, jouini2009queueing, jouini2011call, allon2011impact, allon2011we, ibrahim2015does, whitt1999improving} and references therein analyzes the impact of delay announcements on the queueing process and the abandonment process of the system.  Lastly, the work of \citet{hui1996tell, hul1997impact, pruyn1998effects, munichor2007numbers, sarel1998managing, taylor1994waiting} explores the behavioral aspect of customer waiting and how delays affect customer decisions.  This paper is concerned about the impact of time varying arrival rates and delayed information on the queue length process.  Thus, it is mostly related to the work by \citet{armony2004customer, guo2007analysis, hassin2007information, armony2009impact, guo2009impacts, jouini2009queueing, jouini2011call, allon2011impact, allon2011we, ibrahim2015does, whitt1999improving, armony2015patient, dong2015impact}.  

More recently, there is also research that considers how information can impact the dynamics of queueing systems.  Work by \citet{jennings2015comparisons} compares ticket queues with standard queues.  In a ticket queue, the manager of the queue is unaware of when a customer abandons and is only notified of the abandonment when the customer would have entered service.  This artificially inflates the queue length process and \citet{jennings2015comparisons} determines how much the queue length is inflated because of this loss of information.  However, this work does not consider the aspect of choice and and delays in publishing the information to customers, which is the case in many healthcare and transportation settings.   

This paper analyzes two deterministic queueing models, which describe the dynamics of customer choice and delayed queue length information.  In the first model, the customer receives information about the queue length which is delayed by a parameter $\Delta$.  In the second model, we use a moving average of the queue length over the time interval $\Delta$ to represent the queue length information given to the customer.  The models that we analyze are identical to the models that were analyzed in \citet{pender2016managing}, however, in this paper we add a time varying arrival rate, which is a significant generalization.  This is because queues with time varying arrival rates are notoriously difficult to analyze since many of the standard techniques do not apply.  

However, in this paper, we apply asymptotic analysis techniques like matched asymptotic expansions and the two variable expansion method to analyze our new time varying queueing systems with delayed information.  We show in both models that when the time varying arrival is sinusoidal and the sinusoidal part is small, then the time varying part of the arrival rate does not affect the stability of the queueing dynamics unless the frequency of the arrival rate is twice that of the oscillation frequency.  Our main results in this work represents a novel contribution to the literature in queueing theory and dynamical systems because many real-world queueing systems have time varying rates and it is important to understand when the time varying arrival rate will have an affect on the stability dynamics of the system.  Showing that when the time varying amplitude is small relative to the base arrival rate and the frequency of the arrival rate is not twice that of the critical delay, then the stability dynamics are roughly identical to the non-time varying case.

%**************************************************************************
%**************************************************************************
 \subsection{Main Contributions of Paper}

The contributions of this work can be summarized as follows.    
\begin{itemize}
\item We analyze two two-dimensional fluid models with time varying arrival rates that incorporate customer choice based on delayed queue length information.  In the first model, the information provided to the customer is the queue length delayed by a constant $\Delta$ and in the second model, the information provided to the customer is a moving average of the queue length over the a time window of size $\Delta$.  We show that the impact of the time varying arrival does not shift the value of the critical delay unless the frequency of the arrival rate is twice that of the critical frequency.   
\item We show in both models using the method of multiple time scales that the critical delay, which determines the stability of the delay differential equations, can be shifted by the incorporation of time varying arrival rates.  We also determine the size and impact of this shift.  
\end{itemize} 

%**************************************************************************
%**************************************************************************

\subsection{Organization of Paper}

The remainder of this paper is organized as follows. Section~\ref{sec_CD} gives a brief overview of the infinite server queue with time varying arrival rates and describes the constant delay fluid model.  Using asymptotic expansions, we derive the critical delay threshold under which the queues are balanced if the delay is below the threshold and the queues are asynchronized if the delay is above the threshold.  We also show that the increased or decreased stability because of the time varying arrival rates depends on the sign of the amplitude.  Section~\ref{sec_MA} describes a constant moving average delay fluid model.  Using similar asymptotic expansions, we derive the critical delay threshold under which the queues are balanced if the delay is below the threshold and the queues are asynchronized if the delay is above the threshold in the case of time varying arrival rates.  Finally in Section~\ref{sec_conclusion}, we conclude with directions for future research related to this work.

%**************************************************************************
%**************************************************************************

\section{Constant Delay Fluid Model }\label{sec_CD}

In this section, we present a fluid model with customer choice based on the queue length with a constant delay.  Thus, we begin with two infinite-server queues operating in parallel, where customers choose which queue to join by taking the size of the queue length into account. However, we add the twist that the queue length information that is reported to the customer is delayed by a constant $\Delta$.  Therefore, the queue length that the customer receives is actually the queue length $\Delta$ time units in the past.  Our choice model is identical to that of a Multinomial Logit Model (MNL) where the utility for being served in the $i^{th}$ queue with delayed queue length  $Q_i(t-\Delta)$ is $u_i(Q_i(t-\Delta)) =  − Q_i(t-\Delta)$. Thus, our deterministic queueing model with customer choice, delayed information, and time varying arrival rates can be represented by the two dimensional system of delay differential equations
\begin{eqnarray}
\shortdot{q}_1(t) &=& \lambda(t)  \cdot \frac{\exp(-q_1(t-\Delta))}{\exp(-q_1(t-\Delta)) + \exp(-q_2(t-\Delta))} - \mu q_1(t) \label{ddecd1} \\
\shortdot{q}_2(t) &=& \lambda(t)  \cdot \frac{\exp(-q_2(t-\Delta))}{\exp(-q_1(t-\Delta)) + \exp(-q_2(t-\Delta))} - \mu q_2(t) \label{ddecd2}
\end{eqnarray}
where we assume that $q_1(t)$ and $q_2(t)$ start with different initial functions $\varphi_1(t)$ and $\varphi_2(t)$ on the interval $[-\Delta,0]$, $\lambda(t)$ is the total arrival rate to both queues, and $\mu$ is the service rate of each queue.  

\begin{remark}
When the two delay differential equations are started with the same initial functions, they are identical for all time because of the symmetry of the problem.  Therefore, we will start the system with non-identical initial conditions so the problem is no longer trivial and the dynamics are not identical.  
\end{remark}

In the constant delay model, it is critical to understand the case when the arrival rate is constant and does not depend on time.  In \citet{pender2016managing}, the authors show that the critical delay can be determined from the model parameters and the following theorem is from \citet{pender2016managing}.

\begin{theorem}
For the constant delay choice model, the critical delay parameter is given by the following expression
\begin{equation}
\Delta_{cr}(\lambda, \mu)=\frac{2 \arccos(-2\mu/\lambda)}{ \sqrt{\lambda^2-4\mu^2}}.
\end{equation}

\begin{proof}
See \citet{pender2016managing}.
\end{proof}
\end{theorem}

However, the model of \citet{pender2016managing} neglects to consider the impact of time varying arrival rates.  Time varying arrival rates are important to incorporate into one's model of queues since real customer behavior is dynamic and is not constant over time.   To this end, we will exploit asymptotic analysis and perturbation methods to obtain some insight on the impact of time varying arrival rates.

\subsection{Understanding the $M_t/M/\infty$ Queue }

Before we analyze the queueing model with customer choice it is important to understand the dynamics of the infinite server queue with a time varying arrival rate since it will be essential to our future analysis.  We know from the work of \citet{eick1993mt, eick1993physics}, that the infinite server queue or the $M_t$/G/$\infty$ queue has a Poisson distribution when initialized at zero or with a Poisson distribution with mean rate $q_{\infty}(t)$ where 

\begin{eqnarray}
q_{\infty}(t) &=& E[Q_{\infty}(t) ] \\
&=& \int^{t}_{-\infty} \overline{G}(t-u) \lambda(u) du \\
&=& E\left[\int^{t}_{t-S} \lambda(u) du \right] \\
&=& E[ \lambda( t - S_e ) ] \cdot E[S]
\end{eqnarray}
where $S$ represents a service time with distribution G, $\overline{G} = 1 - G(t) = 
\mathbb{P}( S > t)$, and $S_e$ is a random variable with distribution that follows the 
stationary excess of residual-lifetime cdf $G_e$, defined by
 \begin{eqnarray}
G_e(t) &\equiv&  \mathbb{P}( S_e < t) = \frac{1}{E[S]} \int^{t}_{0} \overline{G}(u) du, \ 
\ \ t \geq 0.
\end{eqnarray}

The exact analysis of the infinite server queue is often useful since it represents the dynamics of the queueing process if there were an unlimited amount of resources to satisfy the demand process.  Moreover, as observed in \citet{pender2014poisson}, when the service time distribution is exponential, the mean of the queue length process $q^{\infty}(t)$ is the solution to the following ordinary linear differential equation
\begin{eqnarray}
\shortdot{q}_{\infty}(t) = \lambda(t) - \mu \cdot q_{\infty}(t).
\end{eqnarray}

\begin{theorem}\label{MtMtinf}
The solution to the mean of the $M_t/M_t/\infty$ queue with initial value $q_0$ is given by 
\begin{eqnarray}
E[Q_{\infty}(t)] &=& q_{\infty}(t) \\
&=& q_0 \cdot \exp\left\{ - \int^{t}_{0} \mu(s) ds \right\}  \\
&&+ \left( \exp\left\{ - \int^{t}_{0} \mu(s) ds \right\}  \cdot \left( \int^{t}_{0} \lambda(s) \exp\left\{ \int^{s}_{0} \mu(r) dr \right\}  ds \right) \right) .
\end{eqnarray}

\begin{proof}
We can exploit the fact that the mean of the Markovian time varying infinite server queue solves a linear ordinary differential equation.  Therefore, we can use standard ode theory to find the mean of the infinite server queue.  For more details see for example, \citet{pender2014poisson}.

\end{proof}
\end{theorem}

\begin{corollary}
In the special case where $q_0 = 0$, $\mu$ is constant, and $\lambda(t) = \lambda + \lambda \cdot \alpha \sin(\gamma t)$, the mean queue has the following representation 
\begin{eqnarray*}
E[Q_{\infty}(t)] &=&  \frac{\lambda}{\mu} \cdot ( 1 -  \exp(-\mu t) ) + \frac{\lambda \cdot \alpha }{\mu^2 + \gamma^2} \cdot \left[   \left( \mu \cdot \sin(\gamma t) - \gamma \cdot \cos(\gamma t)  \right) + \exp(-\mu t) \cdot \gamma  \right].
\end{eqnarray*}

Morever, when $t$ is very large, then we have that 
\begin{eqnarray*}
E[Q_{\infty}(t)] &\approx&  \frac{\lambda}{\mu}  + \frac{\lambda \cdot \alpha }{\mu^2 + \gamma^2} \cdot \left[   \mu \cdot \sin(\gamma t) - \gamma \cdot \cos(\gamma t)   \right].
\end{eqnarray*}
\end{corollary}

%**************************************************************************
%**************************************************************************

\subsection{Constant Delay Model with Time Varying Arrivals}

Although the case where the constant delay queueing model has a constant arrival rate $\lambda$, the extension to more complicated arrival functions such as $\lambda(t) = \lambda + \lambda \cdot \alpha \sin(\gamma t)$ are quite difficult to analyze.  However, we can analyze the system when the time varying arrival rate is close to the the constant rate case using perturbation theory.   Thus, we assume that the queue length equations for the constant delay model satisfy the following delay differential equations  

\begin{eqnarray}
\shortdot{q}_1(t) &=& \left(  \lambda + \lambda \cdot \alpha \cdot \epsilon \sin(\gamma t) \right)  \cdot \frac{\exp(-q_1(t-\Delta))}{\exp(-q_1(t-\Delta)) + \exp(-q_2(t-\Delta))} - \mu q_1(t) \label{ddecd11} \\
\shortdot{q}_2(t) &=& \left(  \lambda + \lambda \cdot \alpha \cdot \epsilon \sin(\gamma t) \right)  \cdot  \frac{\exp(-q_2(t-\Delta))}{\exp(-q_1(t-\Delta)) + \exp(-q_2(t-\Delta))} - \mu q_2(t) \label{ddecd21}
\end{eqnarray}
where we assume that $q_1(t)$ and $q_2(t)$ start with different initial functions $\varphi_1(t)$ and $\varphi_2(t)$ on the interval $[-\Delta,0]$ and we assume that $ 0 \leq \alpha \leq 1$ and $\epsilon \ll 1$.  

In order to begin our analysis of the delay differential equations, we need to understand the case where $\epsilon = 0$.  Fortunately, this analysis has been carried out in \citet{pender2016managing} and we give a brief outline of the analysis for the reader's convenience.  The first step to understanding the case when $\epsilon = 0$ to compute the equilibrium in this case.  
 
In our case, the delay differential equations given in Equations \ref{ddecd11} - \ref{ddecd21} are symmetric.  Moreover, in the case where the delay $\Delta =0$, the two equations converge to the same point since in equilibrium each queue will receive exactly one half of the arrivals and the two service rates are identical.  This is also true in the case where the arrival process contains delays in the queue length since in equilibrium, the delayed queue length is equal to the non-delayed queue length.  Thus, we have in equilibrium that 
\begin{equation}
 q_1(t) = q_2(t) = \frac{q_{\infty}(t)}{2} \quad \mathrm{ as \ } t \to \infty.
\end{equation}
and
\begin{equation}
q_1(t-\Delta) = q_2(t-\Delta) = \frac{q_{\infty}(t-\Delta)}{2} \quad \mathrm{ as \ } t \to \infty.
\end{equation}

Now that we know the equilibrium for Equations \ref{ddecd1} - \ref{ddecd2}, we need to understand the stability of the delay differential equations around the equilibrium.  The first step in doing this is to linearize the non-linear delay differential equations around the equilibrium point.  This can be achieved by setting the queue lengths to 
\begin{eqnarray}
q_1(t) &=& \frac{q_{\infty}(t)}{2} + u(t) \label{sub1}\\
q_2(t) &=& \frac{q_{\infty}(t)}{2} - u(t) \label{sub2}
\end{eqnarray}
where $u(t)$ is a pertubation function about the equilibrium point $\frac{q_{\infty}(t)}{2}$.  By substituting Equations \ref{sub1} - \ref{sub2} into Equations \ref{ddecd1} - \ref{ddecd2} respectively and linearizing around the point $u(t) =0$, we have that the perturbation function solves the following delay differential equation 
\begin{eqnarray}
\shortdot{u}(t) &=& -\frac{\lambda}{2} \cdot u(t-\Delta) - \mu \cdot  u(t) \label{pert}.
\end{eqnarray}

Therefore, it only remains for us to analyze Equation \ref{pert} to understand the stability of the constant delay queueing system.

Now we set $u(t) =\exp(i\omega t)$ in Equation \ref{pert} provides the values for $\omega_{cr}$ and $\Delta_{cr}$:
%\begin{equation}
%\cos\omega\Delta = -2\mu/\lambda,~~~~~~~~~~\sin\omega\Delta = 2\omega/ \lambda
%\label{goo3}
%\end{equation}
\begin{equation}
\omega_{cr}=\frac{1}{2}\sqrt{\lambda^2-4\mu^2}
\label{goo4}
\end{equation}
\begin{equation}
\Delta_{cr}=\frac{2\arccos(-2\mu/\lambda)}{\sqrt{\lambda^2-4\mu^2}}
\label{goo5}
\end{equation}
Note that Equation \ref{pert} possesses the special solution for $\Delta=\Delta_{cr}$:
 \begin{equation}
 u(t) = A \cos \omega_{cr} t + B \sin  \omega_{cr} t   \label{goo6}
 \end{equation}
where $A$ and $B$ are arbitrary constants.\\

\subsection{Asymptotic Expansions for Constant Delay Model}

Now that we understand the case where $\epsilon =0$, it remains for us to understand the general case.  One important observation to make is that in the previous subsection, we did not use the arrival rate in any way.  Therefore, the same analysis can be repeated with the time varying arrival rate with no changes.  Following the same steps as in the case $\epsilon =0$, we arrive at the case where we need to analyze the following delay differential equation
\begin{equation}
\shortdot{z}(t) = -\frac{\lambda}{2}\left(1+ \alpha \cdot \epsilon \cdot \sin\gamma t \right) \cdot z(t-\Delta) -\mu z(t) ,~~~~~~~~~\epsilon \ll 1.
\label{cdddetv}
\end{equation}
However, since the arrival rate is not constant this time, we do not have a simple way to find the stability of the equation.  Therefore, we will exploit the fact that the time varying arrival rate is near the constant arrival rate and use the two variable expansion method or the method of multiple time scales developed by \citet{kevorkian2013perturbation}.  

\begin{theorem}
The only resonant frequency $\gamma$ of the time varying arrival rate function for the first-order two variable expansion is $\gamma=2\omega_{cr}$. For this value of $\gamma$, the change in stability occurs at the value $\Delta_{mod}$ where 
 \begin{equation}
 \Delta_{mod}=\Delta_{cr}-\epsilon\sqrt{\frac{\alpha^2}{\lambda ^2-4 \mu ^2}}.
 \end{equation}
 
 \begin{proof}

We expand time into two variables $\xi$ and $\eta$ that represent regular and slow time respectively i.e.
\begin{equation}
\xi = t \mbox{~~~(regular time) ~~~and~~~~} \eta=\epsilon t \mbox{~~~(slow time)}.
\end{equation}
Therefore, $z(t)$ now becomes $z(\xi,\eta)$, and 
\begin{equation}
\shortdot{z}(t) = \frac{dz}{dt} = \frac{\partial z}{\partial \xi}\frac{d\xi}{dt}+\frac{\partial z}{\partial \eta}\frac{d\eta}{dt}
=\frac{\partial z}{\partial \xi}+\epsilon \frac{\partial z}{\partial \eta}  \label{foo0}
\end{equation}
Moreover, we have that 
 \begin{equation}
 z(t-\Delta) = z(\xi-\Delta,\eta-\epsilon\Delta)
 \label{foo1}
 \end{equation}
In discussing the dynamics of \ref{cdddetv}, we will detune the delay $\Delta$ off of its critical value:
 \begin{equation}
\Delta = \Delta_{cr} + \epsilon \Delta_1 + O(\epsilon^2)
\label{foo2}
 \end{equation}
Substituting Equation \ref{foo2} into Equation \ref{foo1} and expanding term by term, we get
 \begin{equation}
z(t-\Delta)=\bar z-\epsilon\Delta_1\frac{\partial \bar z}{\partial \xi}-\epsilon\Delta_{cr}\frac{\partial \bar z}{\partial \eta}+O(\epsilon^2)
 \label{foo3}
 \end{equation}
where $$\bar z=z(\xi-\Delta_{cr},\eta).$$
Equation \ref{cdddetv} becomes, neglecting terms of $O(\epsilon^2)$,
 \begin{equation}
 \frac{\partial z}{\partial \xi}+\epsilon \frac{\partial z}{\partial \eta}=-\mu z
-\frac{\lambda}{2}\left(1+ \alpha \cdot \epsilon \cdot \sin\gamma t \right)
\left(\bar z-\epsilon\Delta_1\frac{\partial \bar z}{\partial \xi}-\epsilon\Delta_{cr}\frac{\partial \bar z}{\partial \eta}\right)
 \label{foo4}
 \end{equation}
Now we expand $z$ in a power series in $\epsilon$:
 \begin{equation}
z = z_{0} + \epsilon z_{1} + O(\epsilon^2)
 \label{foo5}
 \end{equation}
 
Substituting (\ref{foo5}) into (\ref{foo4}), collecting terms, and equating similar powers of $\epsilon$, we get
\begin{align}
  \frac{\partial z_{0}}{\partial\xi}+\mu z_{0}+ \frac{\lambda}{2} {\bar z_0} &= 0  \label{foo6} \\
  \frac{\partial z_{1}}{\partial\xi}+\mu z_{1}+ \frac{\lambda}{2} {\bar z_1} &= -\frac{\partial z_0}{\partial\eta} +  \frac{\lambda}{2} \left( \Delta_1  \frac{\partial {\bar z_0}}{\partial\xi}
+   \Delta_{cr} \frac{\partial {\bar z_0 }}{\partial\eta} -   \alpha  {\bar z_0} \sin \gamma \xi \right)
  \label{foo7}
\end{align}
  
  Equation \ref{foo6} has the solution given in Equation \ref{goo6}:
 \begin{equation}
 z_0 = A(\eta) \cos \omega_{cr} \xi + B(\eta) \sin  \omega_{cr} \xi .
 \label{foo8}
 \end{equation}
 The functions $A(\eta)$ and $B(\eta)$ give the slow flow of the system. We find differential equations on $A(\eta)$ and $B(\eta)$ by substituting Equation \ref{foo8} into \ref{foo7} and eliminating the resonant terms.
 
The next step is to substitute (\ref{foo8}) into (\ref{foo7}).  The quantity $\bar z_0$ in (\ref{foo7}) may be conveniently computed 
from the following expression, obtained from (\ref{foo6}):
\begin{align}
  \bar z_0 &= \frac{2}{\lambda} \cdot \left(-\mu z_0-\frac{\partial z_0}{\partial\xi}\right) \nonumber \\
  &= \frac{2}{\lambda} \cdot \left[-(\mu A+\omega_{cr}B) \cos \omega_{cr} \xi  +(\omega_{cr} A - \mu B) \sin \omega_{cr} \xi \right] \label{foo9}
\end{align}
 Therefore, we have the following expressions for the terms in Equation \ref{foo7}
     \begin{align}
\frac{\partial z_0}{\partial\eta} &= A' \cdot \cos(\omega_{cr} \xi ) + B' \cdot \sin(\omega_{cr} \xi )  \\
\frac{\partial {\bar z_0}}{\partial\xi}&= \frac{2 \cdot \omega_{cr} }{\lambda} \left[  ( A \omega_{cr} - \mu B )  \cdot \cos(\omega_{cr} \xi ) + (\mu A + B \omega_{cr}) \cdot  \sin(\omega_{cr} \xi )  \right]  \\
\frac{\partial {\bar z_0 }}{\partial\eta} &=  -\frac{2}{\lambda} \left[  ( \mu A' + B' \omega_{cr}) \cdot  \cos(\omega_{cr} \xi ) + (\mu B' - A' \omega_{cr})  \cdot \sin(\omega_{cr} \xi )  \right] \\ 
 \alpha  {\bar z_0}  \sin \gamma \xi  &= -\alpha \cdot \left[  ( \mu A + B \omega_{cr}) \cdot  \cos(\omega_{cr} \xi ) + (\mu B - A \omega_{cr})  \cdot \sin(\omega_{cr} \xi )  \right] \cdot \sin( \gamma \xi) \nonumber \\
 &= \frac{\alpha}{2} (A\omega_{cr}-B\mu)\left[\cos((\gamma-\omega_{cr})\xi) - \cos((\gamma+\omega_{cr})\xi) \right] \nonumber \\
 &\quad - \frac{\alpha}{2}(A\mu+B\omega_{cr})\left[ \sin((\gamma-\omega_{cr})\xi) + \sin((\gamma+\omega_{cr})\xi) \right]
  \label{foo100}
  \end{align}
 
Thus, after substituting Equation \ref{foo8} into \ref{foo7} and applying angle-sum identities, the only terms involving $\gamma$ are of the form 
\begin{equation}
\cos((\gamma\pm\omega_{cr})\xi),\quad \sin((\gamma\pm\omega_{cr})\xi)
\end{equation}
 Notice that $\gamma=2\omega_{cr}$ is the only resonant frequency for the arrival function. For any other value of $\gamma$, the terms involving $\gamma$ at $O(\epsilon)$ are non-resonant, and the first-order two-variable expansion method does not capture any effect from the time-varying arrival function. This 2 to 1 resonance is a similar phenomenon to that arising from ordinary differential equations involving parametric excitation, see for example \citet{ng2003nonlinear, ruelas2012nonlinear}.
 Therefore, we set $\gamma=2\omega_{cr}$, and Equation \ref{foo7} becomes
 
 \begin{align}
   \frac{\partial z_{1}}{\partial\xi}+\mu z_{1}+ \frac{\lambda}{2} {\bar z_1} &= \left[c_1 A'(\eta) + c_2 B'(\eta) + c_3 A(\eta) + c_4 B(\eta)\right]\cos(\omega_{cr}\xi) \nonumber \\
   &\quad + \left[c_5 A'(\eta) + c_6 B'(\eta) + c_7 A(\eta) + c_8 B(\eta)\right]\sin(\omega_{cr}\xi) \nonumber \\
   &\quad +\text{ non-resonant terms} \label{foo10}
 \end{align}
where
\begin{align}
  c_1 &= 1 + \mu\Delta_{cr},  &  c_2 &= \Delta_{cr}\omega_{cr}, &  c_3 &= \frac{\alpha\omega_{cr}}{2}-\Delta_1\omega_{cr}^2, & c_4 &= -\frac{\alpha\mu}{2}+\Delta_1\mu\omega_{cr} \\
  c_5 &= -\Delta_{cr}\omega_{cr}, & c_6 &= 1+\mu\Delta_{cr}, & c_7 &= -\frac{\alpha\mu}{2}-\Delta_1\mu\omega_{cr}, & c_8 &= -\frac{\alpha\omega_{cr}}{2}+\Delta_1\omega_{cr}^2
\end{align}
Elimination of secular terms gives the slow flow:

\begin{eqnarray}
\label{hoo12}
\frac{dA}{d\eta}=K_1 A(\eta) + K_2 B(\eta) \\
\frac{dB}{d\eta}=K_3 A(\eta) + K_4 B(\eta)
\label{hoo13}
\end{eqnarray}

where
\begin{align}
  K_1 &= -\frac{\omega_{cr} (2 \alpha  \Delta_{cr} \mu +\alpha -2 \Delta_1 \omega_{cr})}{2 \left(\Delta_{cr}^2 \omega_{cr}^2+(\Delta_{cr} \mu +1)^2\right)} \label{k1}\\
  K_2 &= \frac{\alpha  \left(\Delta_{cr} \mu ^2-\Delta_{cr} \omega_{cr}^2+\mu \right)-2 \Delta_1 \omega_{cr} \left(\Delta_{cr} \mu ^2+\Delta_{cr} \omega_{cr}^2+\mu \right)}{2 \left(\Delta_{cr}^2 \omega_{cr}^2+(\Delta_{cr} \mu +1)^2\right)} \\
  K_3 &= \frac{\alpha  \left(\Delta_{cr} \mu ^2-\Delta_{cr} \omega_{cr}^2+\mu \right)+2 \Delta_1 \omega_{cr} \left(\Delta_{cr} \mu ^2+\Delta_{cr} \omega_{cr}^2+\mu \right)}{2 \left(\Delta_{cr}^2 \omega_{cr}^2+(\Delta_{cr} \mu +1)^2\right)} \\
  K_4 &= \frac{\omega_{cr} (2 \alpha  \Delta_{cr} \mu +\alpha +2 \Delta_1 \omega_{cr})}{2 \left(\Delta_{cr}^2 \omega_{cr}^2+(\Delta_{cr} \mu +1)^2\right)} \label{k4}
\end{align}
 The equilibrium point $A(\eta)=B(\eta)=0$ of the slow flow corresponds to a periodic solution for $z_0$, and the stability of the equilibrium corresponds to the stability of that periodic solution. The stability is determined by the eigenvalues of the matrix
  \begin{equation}
  K=
 \left[
 \begin{array}{ccc}
 K_1 & K_2\\
 K_3 & K_4
 \end{array}
 \right]
 \label{kmatrix}
 \end{equation}
 If both eigenvalues have negative real part, the equilibrium is stable.
 Since the eigenvalues are cumbersome to work with directly, we use the Routh-Hurwitz stability criterion: 
 
 Denote the characteristic polynomial of $K$ by 
 \begin{equation}
 \det(K-r I) = a_0 + a_1 r + a_2 r^2 = 0 \label{charpoly}
 \end{equation}
 Then both eigenvalues have negative real part if and only if all the coefficients satisfy $a_i > 0$. From Equations \ref{k1}-\ref{k4} and \ref{charpoly}, and using the expression for $\omega_{cr}$ from Equation \ref{goo4},  we have
\begin{align}
  a_0 &= \frac{\left(\mu ^2+\omega_{cr}^2\right) \left(4 \Delta_1^2 \omega_{cr}^2-\alpha ^2\right)}{4 \left(\Delta_{cr}^2 \omega_{cr}^2+(\Delta_{cr}\mu +1)^2\right)} 
  = \frac{\lambda ^2 \left(\Delta_1^2 \left(\lambda ^2-4 \mu ^2\right)-\alpha ^2\right)}{4 \left(\Delta_{cr}^2 \lambda ^2+8 \Delta_{cr} \mu +4\right)} \label{a0val}\\
  a_1 &= -\frac{2 \Delta_1 \omega_{cr}^2}{\Delta_{cr}^2 \omega_{cr}^2+(\Delta_{cr} \mu +1)^2} 
  =  -\frac{2 \Delta_1 \left(\lambda ^2-4 \mu ^2\right)}{\Delta_{cr}^2 \lambda ^2+8 \Delta_{cr}\mu +4} \label{a1val} \\
  a_2 &= 1 \label{a2val}
\end{align}
 Recall that $\omega_{cr}$ is real and positive only if $\lambda>2\mu$. So, using this restriction, we find that all of the $a_i$ are positive if and only if
 \begin{equation}
 \Delta_1 < -\frac{|\alpha|}{2\omega_{cr}} = -\sqrt{\frac{\alpha^2}{\lambda ^2-4 \mu ^2}}. \label{delta1}
 \end{equation}
 
 Note that we can recover the case with no resonant forcing by setting the forcing amplitude $\alpha=0$. With no forcing, the periodic solution for $z_0$ becomes unstable at $\Delta=\Delta_{cr}$, but with resonant forcing, the change of stability occurs when
 \begin{equation}
 \Delta_{mod}=\Delta_{cr}-\epsilon\sqrt{\frac{\alpha^2}{\lambda ^2-4 \mu ^2}}.
 \end{equation}

 \end{proof}
\end{theorem}

\subsection{Numerics for Constant Delay Queueing Model}
\begin{figure}
\captionsetup{justification=centering}
		\hspace{-.35in}~\includegraphics[scale = .22]{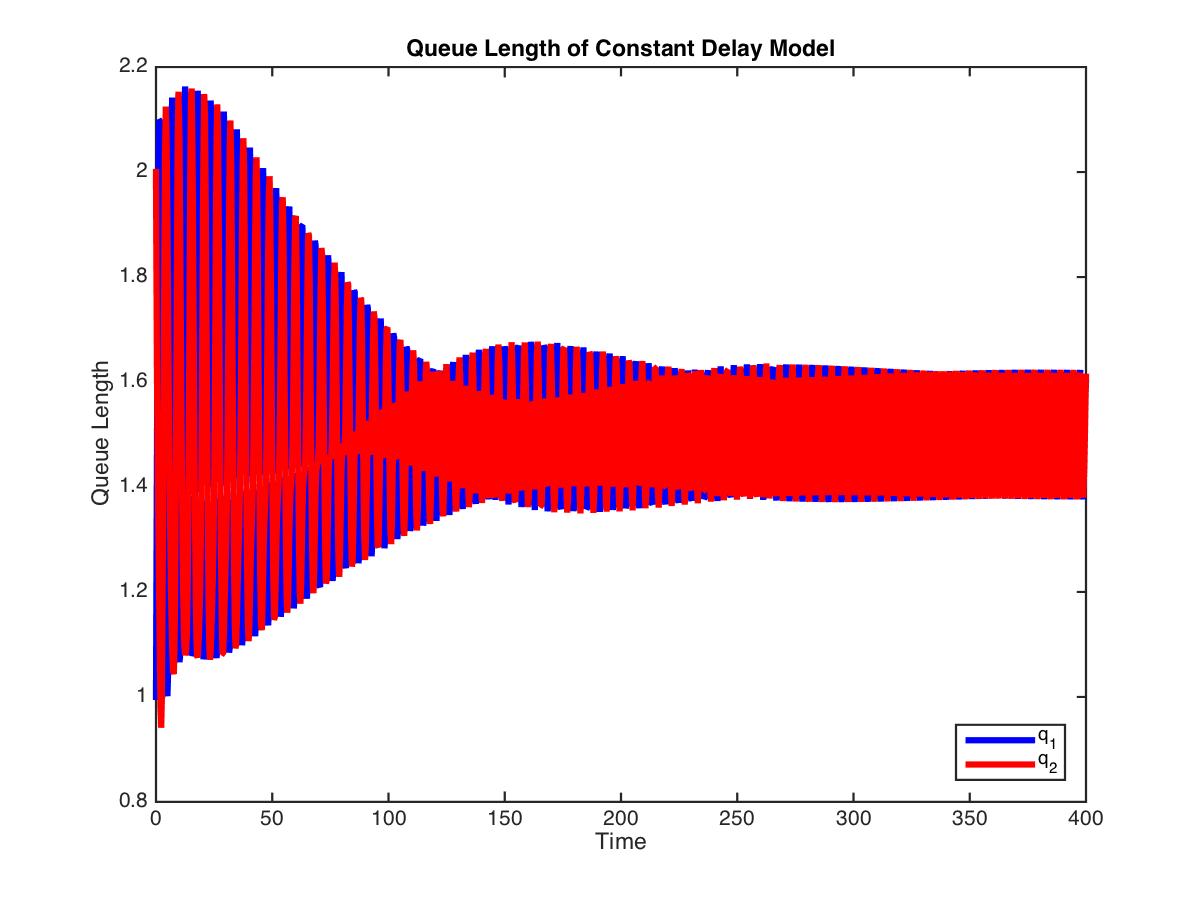}~\hspace{-.4in}~\includegraphics[scale = .22]{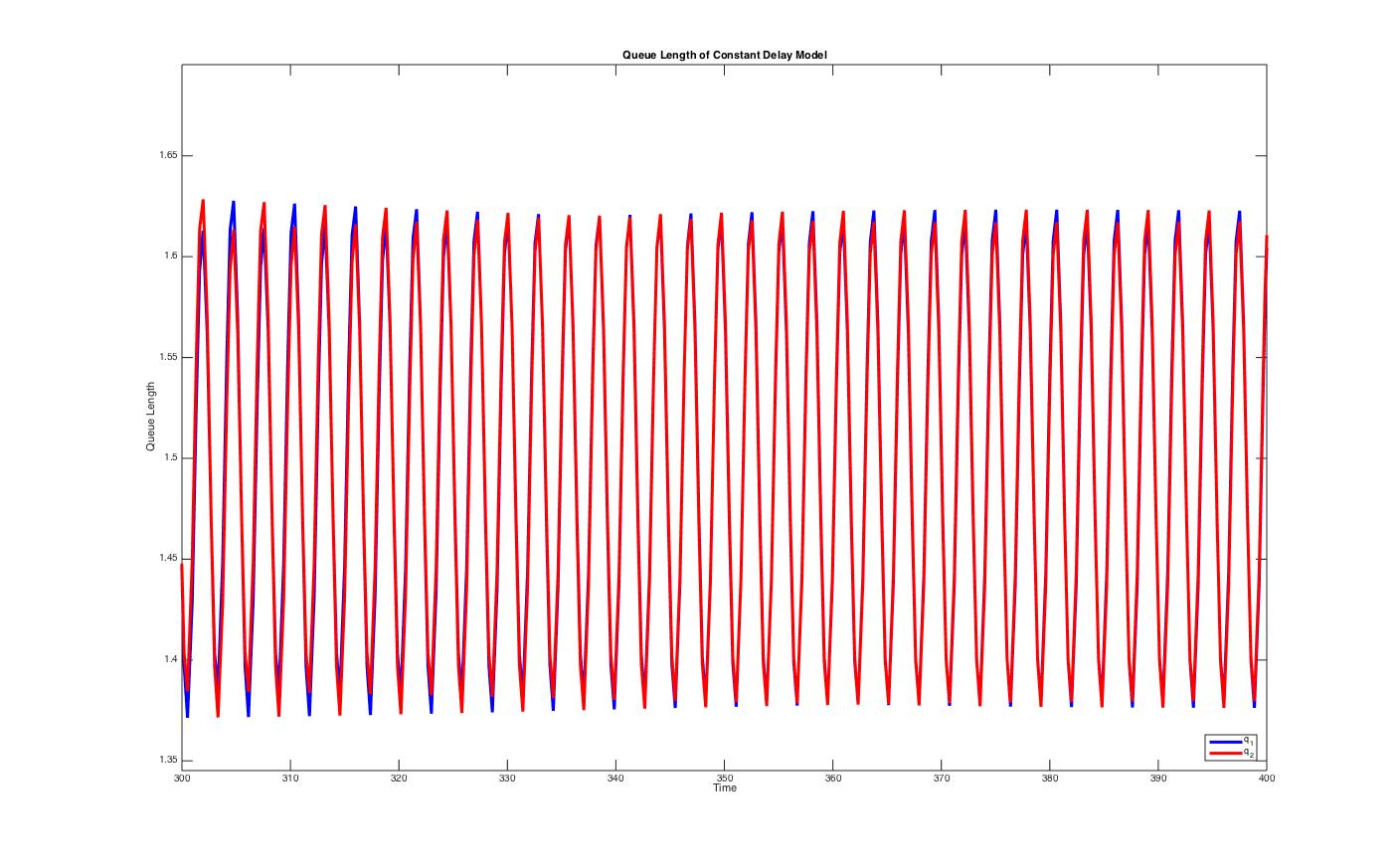}
\caption{ $\Delta_{cr}$ = 2.0577, $\Delta_{mod}$  = 1.9682. \\
$\lambda = 3$, $\mu =1$, $\alpha =1$, $\epsilon = .2$, $\gamma = \sqrt{5}$, $\Delta$ = 1.947, $\varphi_1([-\Delta,0])=1$,$\varphi_2([-\Delta,0])=2$ } \label{Fig5}
\end{figure}

\begin{figure}
\captionsetup{justification=centering}
		\hspace{-.35in}~\includegraphics[scale = .22]{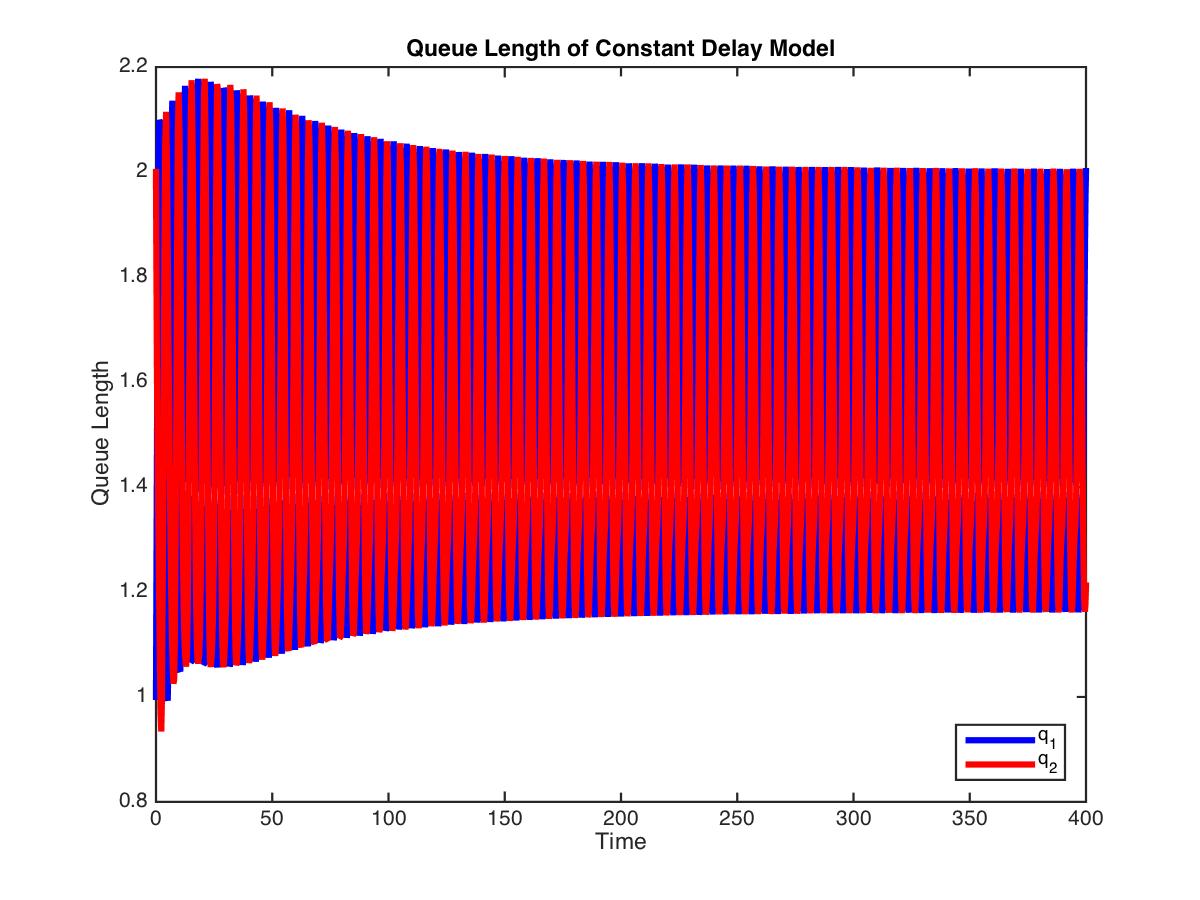}~\hspace{-.3in}~\includegraphics[scale = .22]{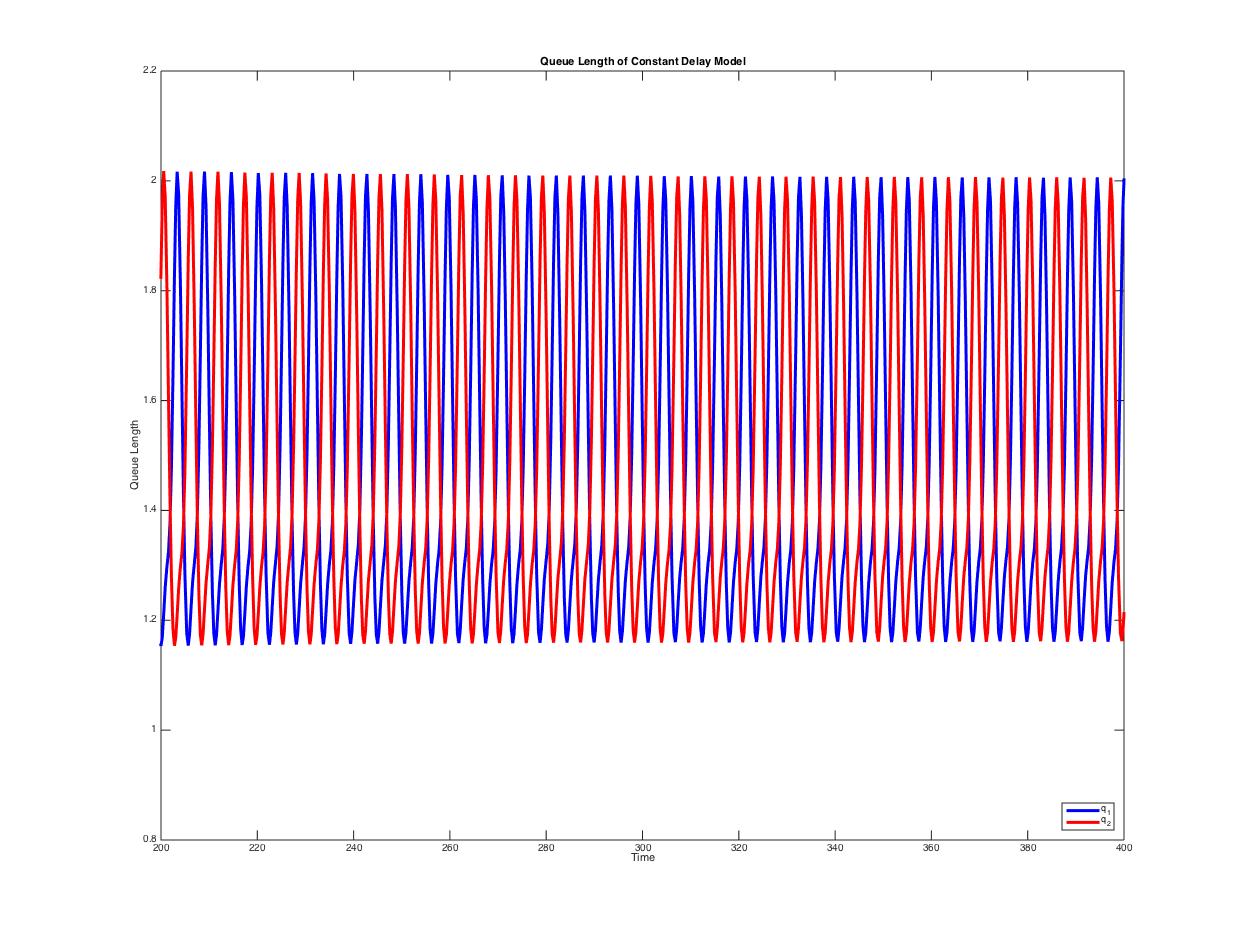}
\caption{ $\Delta_{cr}$ = 2.0577, $\Delta_{mod}$  = 1.9682. \\
$\lambda = 3$, $\mu =1$, $\alpha =1$, $\epsilon = .2$, $\gamma = \sqrt{5}$, $\Delta$ = 1.977, $\varphi_1([-\Delta,0])=1$,$\varphi_2([-\Delta,0])=2$} \label{Fig6}
\end{figure}

 In this section, we numerically integrate the delay two examples of delay differential equations with costant delays and compare the asymptotic results for determining the Hopf bifurcation that occurs.  On the left of Figure \ref{Fig5} we numerically integrate the two queues and plot the queue length as a function of time.  In this example our lag in information is given by $\Delta =1.947$.  We see that the two equations eventually converge to the same limit as time is increased towards infinity.  This implies that the system is stable and no oscillations or asynchrous dynamics will occur due to instability in this case.  On the right of Figure \ref{Fig5} is a zoomed in version of the figure on the left.   It is clear that the two delay equations are converging towards one another and this system is stable.  However, in Figure \ref{Fig6}. we use the same parameters, but we make the lag in information $\Delta = 1.977$.  This is below the critcal delay in the constant case and above the modified critical delay when the time varying arrival rate is taken into account.  On the right of Figure \ref{Fig6}, we display a zoomed in version of the figure on the left. We see that in this case the two queues oscillate and asynchronous behavior is observed.  Thus, the asymptotic analysis performed works well at predicting the change in stability.    

\begin{figure}
\captionsetup{justification=centering}
		\hspace{-.35in}~\includegraphics[scale = .22]{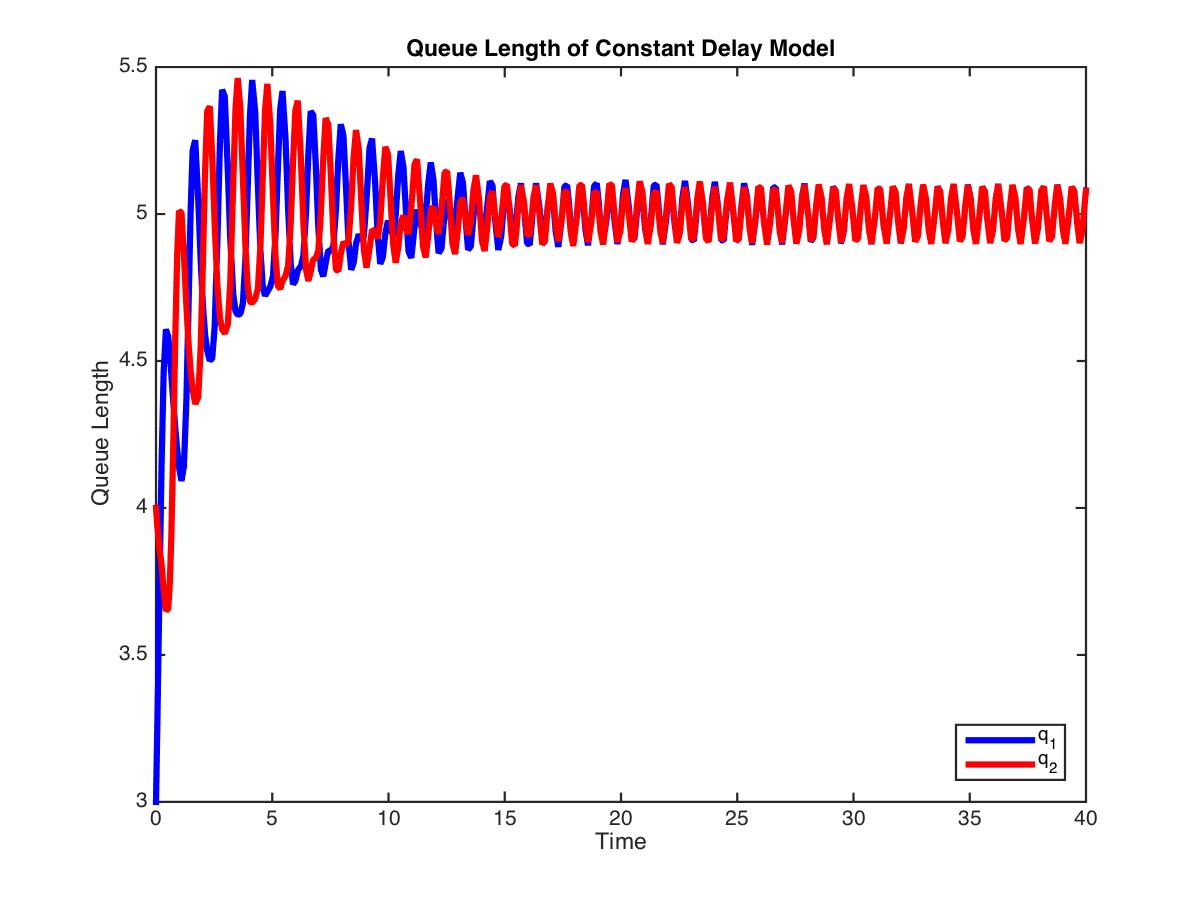}~\hspace{-.4in}~\includegraphics[scale = .22]{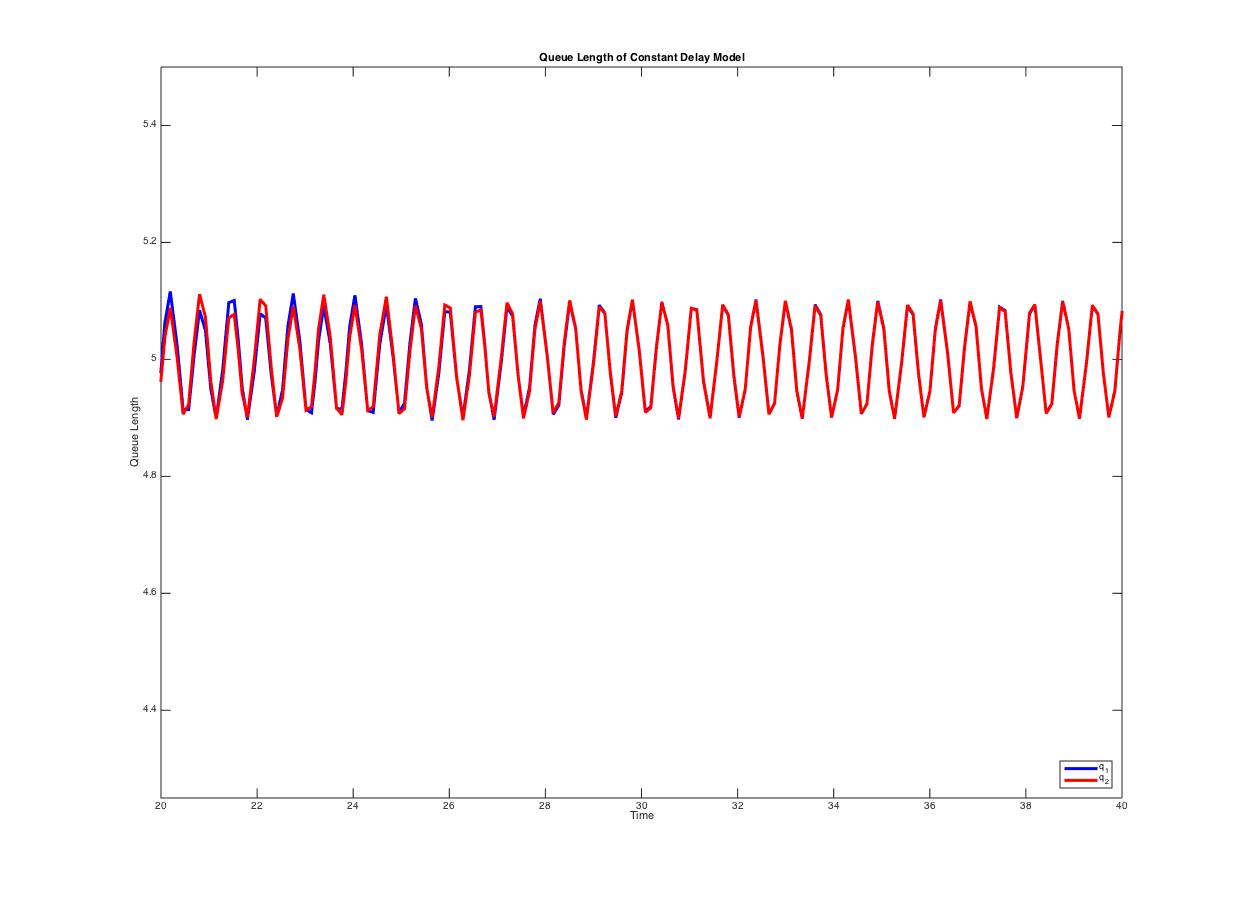}
\caption{ $\Delta_{cr}$ = .3617, $\Delta_{mod}$  = .3413. \\
$\lambda = 10$, $\mu =1$, $\alpha =1$, $\epsilon = .2$, $\gamma = \sqrt{96}$, $\Delta$ = .33, $\varphi_1([-\Delta,0])=3$,$\varphi_2([-\Delta,0])=4$ } \label{Fig7}
\end{figure}

\begin{figure}
\captionsetup{justification=centering}
		\hspace{-.35in}~\includegraphics[scale = .22]{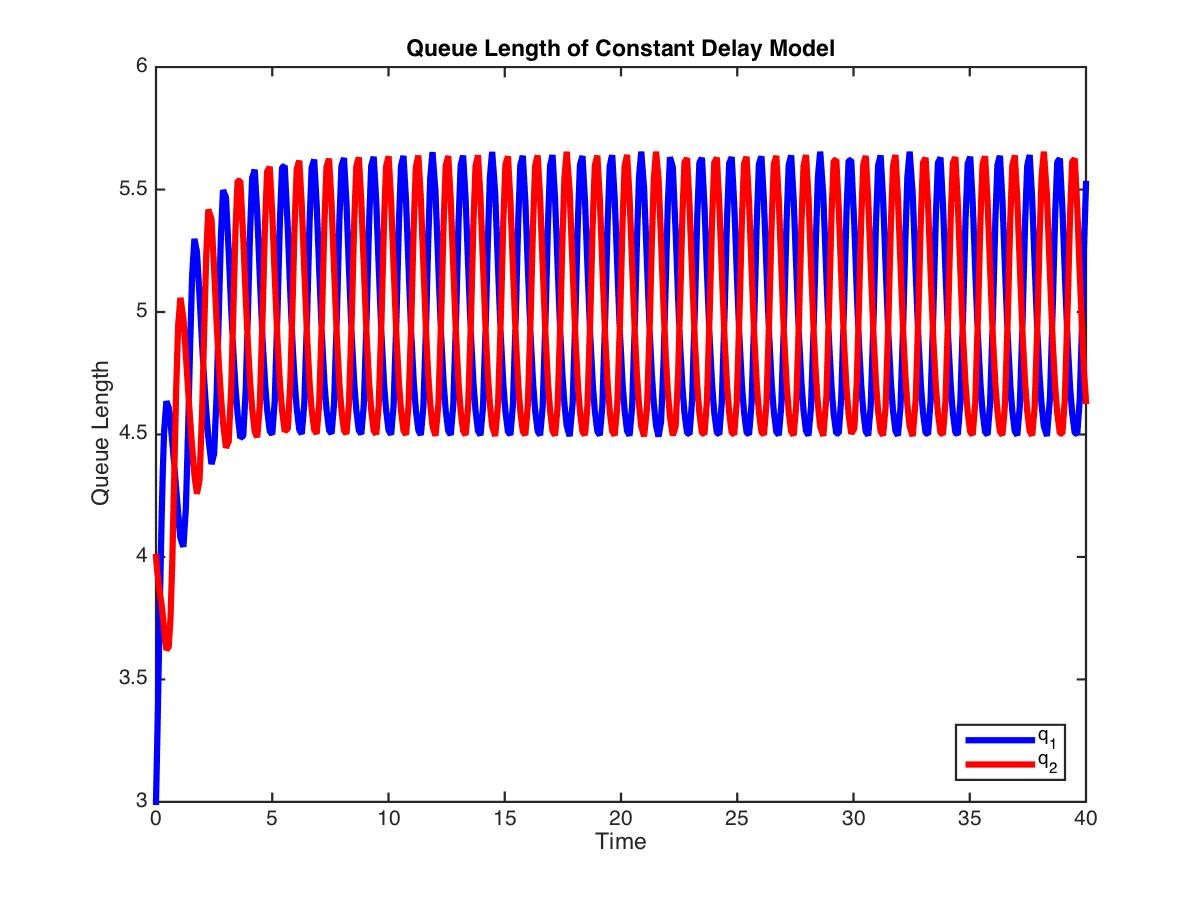}~\hspace{-.3in}~\includegraphics[scale = .22]{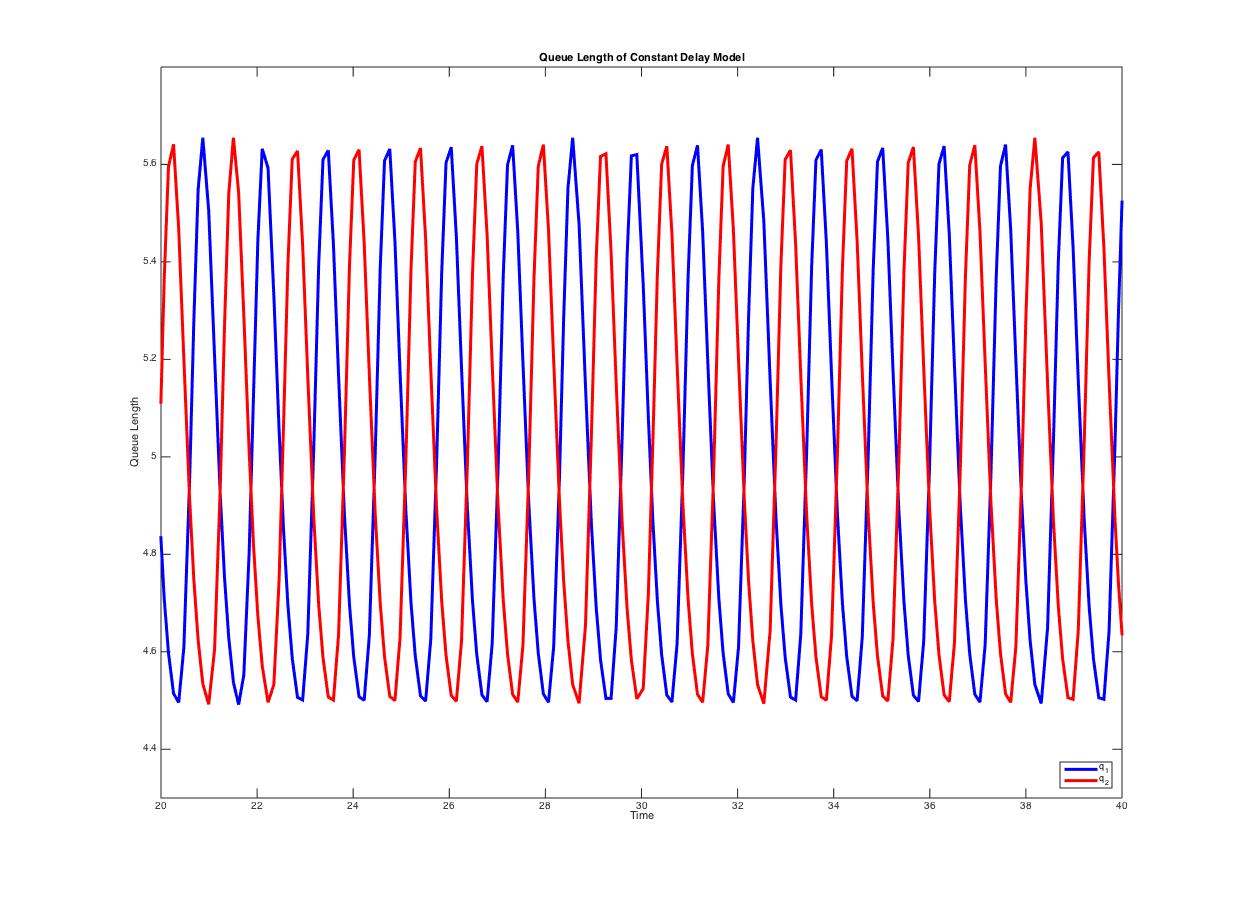}
\caption{ $\Delta_{cr}$ = .3617, $\Delta_{mod}$  = .3413. \\
$\lambda = 10$, $\mu =1$, $\alpha =1$, $\epsilon = .2$, $\gamma = \sqrt{96}$, $\Delta$ = .35, $\varphi_1([-\Delta,0])=3$,$\varphi_2([-\Delta,0])=4$} \label{Fig8}
\end{figure}

As an additional numerical example, on the left of Figure \ref{Fig7}, we numerically integrate the two queues and plot the queue length as a function of time.  In this example our lag in information is given by $\Delta =.33$.  We see that the two equations eventually converge to the same limit as time is increased towards infinity.  This implies that the system is stable and no oscillations or asynchrous dynamics will occur due to instability in this case.  On the right of Figure \ref{Fig7} is a zoomed in version of the figure on the left.   It is clear that the two delay equations are converging towards one another and this system is stable.  However, in Figure \ref{Fig8}. we use the same parameters, but we make the lag in information $\Delta = .35$.  This is below the critical delay in the constant case and above the modified critical delay when the time varying arrival rate is taken into account.  On the right of Figure \ref{Fig8}, we display a zoomed in version of the figure on the left. We see that in this case the two queues oscillate and asynchronous behavior is observed.  Thus, the asymptotic analysis performed works well at predicting the change in stability.

%
%\begin{figure}
%\captionsetup{justification=centering}
%\begin{center}
%		\hspace{-.35in}~\includegraphics[scale = .22]{Time_Varying_Jamol_1.jpg}~\hspace{-.3in}~\includegraphics[scale = .22]{Time_Varying_Jamol_2.jpg}\\
%		\hspace{-.35in}~\includegraphics[scale = .22]{Time_Varying_Jamol_3.jpg}~\hspace{-.3in}~\includegraphics[scale = .22]{Time_Varying_Jamol_4.jpg}\\
%		\hspace{-.35in}~\includegraphics[scale = .22]{Time_Varying_Jamol_5.jpg}~\hspace{-.3in}~\includegraphics[scale = .22]{Time_Varying_Jamol_6.jpg}
%		\end{center}
%\caption{Various Plots as a function of $\alpha$.} \label{Fig5}
%\end{figure}

  %**************************************************************************
  %**************************************************************************

\section{Moving Average Delay Fluid Model} \label{sec_MA}

\begin{figure}
	\centering
		\includegraphics[scale=.95]{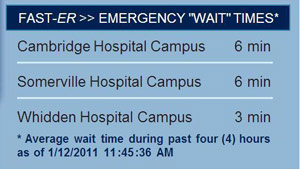}
\caption{Emergency Room Wait Times Via Moving Averages.} \label{FigJFK}
\end{figure}

In this section, we present another fluid model with customer choice and where the delay information presented to the customer is a moving average. This model assumes that customers are informed about the queue length, but in the form of a moving average of the queue length between the current time and $\Delta$ time units in the past.  These types of moving average models are currently used in many healthcare settings such as the one in Figure \ref{FigJFK}.  In Figure \ref{FigJFK}, it is clear that the time information is given in past and is averages over a 4 hour window.  This is partially because patients in healthcare are quite heterogeneous and require different services and attention.  Moreover, the system is not necessary FIFO or FCFS since patients have different priority levels.  Thus, the moving average waiting time indicator might be attractive for these reasons.  
Like in the previous model with constant delays, customers in the moving average model also have the choice to join two parallel infinite server queues and they join according to the same multinomial logit model.   Once again, the extension to more complicated arrival functions such as $\lambda(t) = \lambda + \lambda \cdot \alpha \sin(\gamma t)$ is quite difficult.  However, like in the constant delay setting, we can analyze the system when the time varying arrival rate is close to the the constant rate case using perturbation theory and asymptotics.   Thus, we assume that the queue length equations for the constant delay model satisfy the following delay differential equations  
\begin{equation}
 \lambda(t) \cdot \frac{\exp\left(- \frac{1}{\Delta} \int^{t}_{t-\Delta} q_1(s) ds \right)}{\exp\left(- \frac{1}{\Delta} \int^{t}_{t-\Delta} q_1(s) ds \right) + \exp\left(- \frac{1}{\Delta} \int^{t}_{t-\Delta} q_2(s) ds \right)}
\end{equation}
and join the second queue at rate
\begin{equation}
 \lambda(t)   \cdot \frac{\exp\left(- \frac{1}{\Delta} \int^{t}_{t-\Delta} q_2(s) ds \right)}{\exp\left(- \frac{1}{\Delta} \int^{t}_{t-\Delta} q_1(s) ds \right) + \exp\left( - \frac{1}{\Delta} \int^{t}_{t-\Delta} q_2(s) ds \right)} .
\end{equation} 
Thus, our model for customer choice with delayed information in the form of a moving average can be represented by a two dimensional system of functional differential equations
\begin{align}
\shortdot{q}_1(t) &=  \left( \lambda + \lambda \alpha \epsilon \sin(\gamma t) \right)  \frac{\exp\left(- \frac{1}{\Delta} \int^{t}_{t-\Delta} q_1(s) ds \right)}{\exp\left(- \frac{1}{\Delta} \int^{t}_{t-\Delta} q_1(s) ds \right) + \exp\left(- \frac{1}{\Delta} \int^{t}_{t-\Delta} q_2(s) ds \right)} - \mu q_1(t) \\
\shortdot{q}_2(t) &=  \left( \lambda + \lambda \alpha \epsilon \sin(\gamma t) \right)  \frac{\exp\left(- \frac{1}{\Delta} \int^{t}_{t-\Delta} q_2(s) ds \right)}{\exp\left(- \frac{1}{\Delta} \int^{t}_{t-\Delta} q_1(s) ds \right) + \exp\left(- \frac{1}{\Delta} \int^{t}_{t-\Delta} q_2(s) ds \right)}- \mu q_2(t)
\end{align}
where we assume that $q_1$ and $q_2$ start at different initial functions $\varphi_1(t)$ and $\varphi_2(t)$ on the interval $[-\Delta,0]$.  

On the onset this problem is seemingly more difficult than the constant delay setting since the ratio now depends on a moving average of the queue length during a delay period $\Delta$.  To simplify the notation, we find it useful to define the moving average of the $i^{th}$ queue over the time interval $[t - \Delta, t]$ as
\begin{eqnarray}\label{mai}
 m_i(t,\Delta) = \frac{1}{\Delta} \int^{t}_{t-\Delta} q_i(s) ds.
\end{eqnarray}

This representation of the moving average leads to a key observation where we discover that the moving average itself solves a linear delay differential equation.  In fact, by differentiating Equation \ref{mai} with respect to time, it can be shown that the moving average of the $i^{th}$ queue is the solution to the following delay differential equation
\begin{eqnarray} 
\updot{m}_i(t,\Delta) =   \frac{1}{\Delta} \cdot \left(  q_i(t) - q_i(t - \Delta) \right), \quad i \in \{1,2\}.
\end{eqnarray}

Leveraging the above delay equation for the moving average, we can describe our moving average fluid model with the following four dimensional system of delay differential equations

\begin{eqnarray}
\shortdot{q}_1 &=& \left( \lambda + \lambda \cdot \alpha \cdot \epsilon \cdot \sin(\gamma t) \right) \cdot \frac{\exp\left(- m_1(t) \right) }{\exp\left(- m_1(t) \right) + \exp\left(- m_2(t) \right)} - \mu q_1(t)    \label{ma1}  \\
\shortdot{q}_2 &=& \left( \lambda + \lambda \cdot \alpha \cdot \epsilon \cdot \sin(\gamma t) \right) \cdot \frac{\exp\left(- m_2(t) \right) }{\exp\left(-m_1(t) \right) + \exp\left(- m_2(t) \right)} - \mu q_2(t) \\
\updot{m}_1 &=&  \frac{1}{\Delta} \cdot \left( q_1(t) - q_1(t - \Delta) \right) \\
\updot{m}_2 &=&  \frac{1}{\Delta} \cdot \left( q_2(t) - q_2(t - \Delta) \right) \label{ma2}.
\end{eqnarray}

In the moving average model, it is also critical to understand the case when the arrival rate is constant and does not depend on time.  In \citet{pender2016managing}, the authors show that the critical delay for the moving average model can be determined from the model parameters and the following theorem is from \citet{pender2016managing}.
 
\begin{theorem}
For the moving average fluid model given by Equations \ref{ma1} - \ref{ma2}, the critical delay parameter is the solution to the following transcendental equation
\begin{equation}
\sin\left( \Delta \cdot \sqrt{\frac{\lambda}{\Delta} - \mu^2} \right) + \frac{2 \mu \Delta}{\lambda} \cdot \sqrt{\frac{\lambda}{\Delta} - \mu^2} = 0 .
\end{equation}

\begin{proof}
See \citet{pender2016managing}.
\end{proof}
\end{theorem}

In order to begin our analysis of the delay differential equations with a time varying rate, we need to first understand the case where $\epsilon = 0$ and the arrival rate is constant.  Also like in the constant delay setting, this analysis has been carried out in \citet{pender2016managing} and we give a brief outline of the analysis for the reader's convenience.

The first step to understanding the case when $\epsilon = 0$ to compute the equilibrium in this case.  The first part of the proof is to compute an equilibrium for the solution to the delay differential equations.  In our case, the delay differential equations given in Equations \ref{ma1} - \ref{ma2} are symmetric.  Moreover, in the case where there is no delay, the two equations converge to the same point since in equilibrium each queue will receive exactly one half of the arrivals and the two service rates are identical.  This is also true in the case where the arrival process contains delays in the queue length since in equilibrium, the delayed queue length is equal to the non-delayed queue length.  It can be shown that there is only one equilibrium where all of the states are equal to each other.  One can prove this by substituting $q_2 = \lambda/\mu - q_1$ in the steady state verison of Equation \ref{ma1} and solving for $q_1$.  One eventually sees that $q_1 = q_2$ is the only solution since any other solution does not obey Equation \ref{ma1}.  Thus, we have in equilibrium that 
\begin{eqnarray}
q_1(t) = q_2(t) =  \frac{q_{\infty}(t)}{2 } \quad \mathrm{ as \ } t \to \infty
\end{eqnarray}
and 
\begin{eqnarray}
m_1(t) = m_2(t) = \frac{1}{\Delta} \int^{t}_{t-\Delta} \frac{q_{\infty}(s)}{2} ds \quad \mathrm{ as \ } t \to \infty.
\end{eqnarray}
Now that we know the equilibrium for Equations \ref{ma1} - \ref{ma2}, we need to understand the stability of the delay differential equations around the equilibrium.  The first step in doing this is to set each of the queue lengths to the equilibrium values plus a perturbation.  Thus, we set each of the queue lengths to 
\begin{eqnarray}
q_1(t) &=& \frac{q_{\infty}(t)}{2} + u(t) \label{sub11} \\
q_2(t) &=& \frac{q_{\infty}(t)}{2} - u(t) \\
m_1(t) &=& \frac{1}{\Delta} \int^{t}_{t-\Delta} \frac{q_{\infty}(s)}{2} ds + w(t) \\
m_2(t) &=& \frac{1}{\Delta}\int^{t}_{t-\Delta} \frac{q_{\infty}(s)}{2}ds  - w(t) \label{sub22}
\end{eqnarray}
Substitute Equations  \ref{sub11} - \ref{sub22} into Equations  \ref{ma1} - \ref{ma2} and solve for $\shortdot{q}_\infty$, $\shortdot{u}$ and $\shortdot{w}$. 
\begin{align}
    \shortdot{ q}_\infty &= \lambda+\lambda\alpha\epsilon\sin(\gamma t) - \mu q_\infty (t) \label{qinf-dot} \\
  \shortdot{ u} &= - \frac{\lambda}{2}\left( 1+\alpha\epsilon\sin(\gamma t) \right)\tanh(w(t)) - \mu u(t) \label{u-dot}\\
    \shortdot{ w}  &= \frac{1}{\Delta}\left( u(t) - u(t-\Delta) \right) \label{w-dot}
\end{align}
Equation \ref{qinf-dot} can be solved explicitly, to give the steady-state solution
\begin{equation}\label{qinf-sol}
q_\infty (t) = c e^{-\mu t}+\frac{\lambda}{2} \left(\frac{1}{\mu }+\frac{\alpha  \epsilon  (\mu  \sin  (\gamma  t)-\gamma  \cos  (\gamma  t))}{\gamma ^2+\mu ^2}\right)
\end{equation}
where 
\begin{equation}
c = q_\infty(0) - \frac{\lambda}{2} \left(\frac{1 }{\mu }-\frac{\alpha  \gamma    \epsilon }{\gamma ^2+\mu ^2}\right)
\end{equation}
To determine the stability of the system, we linearize about the point  $u(t) = w(t) = 0$, giving
\begin{align}
    \shortdot{u} &= - \frac{\lambda}{2}\left( 1+\alpha\epsilon\sin(\gamma t) \right)w(t) - \mu u(t) \label{u-lin} \\
    \shortdot{w} &= \frac{1}{\Delta}\left( u(t) - u(t-\Delta) \right) \label{w-lin}
\end{align}

First consider the unperturbed case ($\epsilon=0$):
\begin{align}
    \shortdot{u} &= - \frac{\lambda}{2}w(t) - \mu u(t) \label{u-lin0} \\
    \shortdot{w} &= \frac{1}{\Delta}\left( u(t) - u(t-\Delta) \right) \label{w-lin0}
\end{align}
To study Equations \ref{u-lin0} and \ref{w-lin0}, we set
\begin{eqnarray}
u &=& A \exp(rt) \\
w &=& B \exp(rt) .
\end{eqnarray}
These solutions imply the following relationships between the constants A,B, and r.

\begin{eqnarray}
A r  &=& - \frac{\lambda}{2} B - \mu A \\
B r &=& \frac{1}{\Delta} ( A - A \exp( -r \Delta))
\end{eqnarray}

solving for A yields

\begin{eqnarray}
A   &=& - \frac{\lambda}{2(\mu + r )} B 
\end{eqnarray}
and rearranging yields the following equation for $r$
\begin{eqnarray} \label{r-trans2}
r &=& \frac{\lambda}{2 \Delta \cdot r} ( \exp( -r \Delta) - 1) - \mu .
\end{eqnarray}
Now it remains for us to understand the transition between stable and unstable solutions once again.

To find the transition between stable and unstable solutions, set $r=i\omega$, giving us the following equation
\begin{equation}
i\omega = \frac{\lambda}{2 \Delta i \omega}( \exp(-i \omega \Delta) -1)-\mu.
\end{equation}
Multiplying both sides by $i \omega$ and using Euler's identity, we have that 
\begin{equation} \label{maeqniw}
 \frac{\lambda}{2 \Delta }( \cos( \omega \Delta) - i \sin(\omega \Delta) -1)-\mu i \omega + \omega^2 = 0.
\end{equation}
Writing the real and imaginary parts of Equation \ref{maeqniw}, we get:
\begin{equation}\label{coseqnma}
\cos(\omega \Delta) = 1 - \frac{2 \Delta \omega^2}{\lambda}
\end{equation}
for the real part and 
\begin{equation}\label{sineqnma}
\sin(\omega \Delta) = - \frac{2 \Delta \mu \omega}{\lambda}
\end{equation}
Once again by squaring and adding $\sin\omega \Delta$ and $\cos\omega \Delta$  together, we get:
\begin{equation} \label{omegadef}
\omega = \sqrt{\frac{\lambda}{\Delta} - \mu^2}
\end{equation}

Finally by substituting the expression for $\omega$ into Equations \ref{sineqnma} and \ref{coseqnma} gives us the final expression for the critical delay $\Delta_{cr}$, which is the simultaneous solution to the following transcendental equations:
\begin{equation}\label{crit-2}
   \sin\left( \Delta_{cr}  \sqrt{\frac{\lambda}{\Delta_{cr}} - \mu^2} \right) + \frac{2 \mu \Delta_{cr}}{\lambda}  \sqrt{\frac{\lambda}{\Delta_{cr}} - \mu^2} =0 
\end{equation}
\begin{equation}\label{crit-2a}
   \cos\left( \Delta_{cr}  \sqrt{\frac{\lambda}{\Delta_{cr}} - \mu^2} \right) + 1 - \frac{2\mu^2\Delta_{cr}}{\lambda} =0
\end{equation}

Squaring Equations \ref{crit-2} and \ref{crit-2a} and adding them, we see that they are satisfied simultaneously when
\begin{equation}\label{crit-simul}
2+\left(2-\frac{4 \Delta_{cr} \mu^2}{\lambda }\right) \cos \left(\Delta_{cr} \sqrt{\frac{\lambda }{\Delta_{cr}}-\mu ^2}\right)
+ \frac{4 \Delta_{cr} \mu}{\lambda }\sqrt{\frac{\lambda }{\Delta_{cr}}-\mu ^2} \sin \left(\Delta_{cr} \sqrt{\frac{\lambda }{\Delta_{cr}}-\mu ^2}\right)=0
\end{equation}

%**************************************************************************
%**************************************************************************

\subsection{Asymptotic Expansions for Moving Average Model}

Now that we understand the case where $\epsilon =0$, it remains for us to understand the general case. Recall that we are analyzing the stability of the
linearized system 

\begin{align}
    \shortdot{u} &= - \frac{\lambda}{2}\left( 1+\alpha\epsilon\sin(\gamma t) \right)w(t) - \mu u(t) \tag{\ref{u-lin}} \\
    \shortdot{w} &= \frac{1}{\Delta}\left( u(t) - u(t-\Delta) \right) \tag{\ref{w-lin}}
\end{align}
It is useful to convert the system of two first-order equations to a single second-order equation, by differentiating Equation \ref{u-lin} and substituting in expressions for $w(t)$ and $\shortdot{w}(t)$ from Equations \ref{u-lin} and \ref{w-lin}. We obtain
\begin{align}
\dshortdot{u} = &\left(\frac{\alpha  \gamma  \epsilon  \cos (\gamma  t)}{\alpha  \epsilon  \sin (\gamma  t)+1}-\mu\right)\shortdot{u} + \left( \frac{\alpha  \gamma  \mu  \epsilon  \cos (\gamma  t)}{\alpha  \epsilon  \sin (\gamma  t)+1} -\frac{\lambda + \alpha  \lambda  \epsilon  \sin (\gamma  t)}{2 \Delta } \right)u  \nonumber\\
&+ \left( \frac{\lambda +\alpha  \lambda  \epsilon  \sin (\gamma  t)}{2 \Delta } \right)u(t-\Delta) \label{ddot-u}
\end{align}

However, since the arrival rate is not constant this time, we do not have a simple way to find the stability of the equation.  Therefore, we will exploit the fact that the time varying arrival rate is near the constant arrival rate and use the two variable expansion method.

\begin{theorem}
The only resonant frequency $\gamma$ of the time varying arrival rate function for the first-order two variable expansion is $\gamma=2\omega_{cr}$. For this value of $\gamma$, the change in stability occurs at $\Delta)_{mod}$ where
\begin{equation}
\Delta_{mod}  = \Delta_{cr} \pm \epsilon \sqrt{\frac{\alpha ^2 \Delta _{cr}^2}{\Delta _{cr} \lambda +4 \Delta _{cr} \mu +4}}
\end{equation}
where the sign of the $\epsilon$ term is positive if $\Delta_{cr}>\frac{\lambda -2 \mu }{2 \mu ^2}$ and negative if $\Delta_{cr}<\frac{\lambda -2 \mu }{2 \mu ^2}$.

\begin{proof}

We expand time into two variables $\xi$ and $\eta$ that represent regular and slow time respectively i.e.
\begin{equation}
\xi = t \mbox{~~~(regular time) ~~~and~~~~} \eta=\epsilon t \mbox{~~~(slow time)}.
\end{equation}
Therefore, $u(t)$ now becomes $u(\xi,\eta)$.  Moreover, 
\begin{align}
\shortdot{u} = \frac{du}{dt} &= \frac{\partial u}{\partial\xi}\frac{d\xi}{dt} + \frac{\partial u}{\partial\eta}\frac{d\eta}{dt} \nonumber \\
&= \frac{\partial u}{\partial\xi} + \epsilon\frac{\partial u}{\partial\eta} \label{udot2} \\
\dshortdot{u} = \frac{d^2 u}{dt^2} &= \frac{d}{dt}\left( \frac{\partial u}{\partial\xi} + \epsilon\frac{\partial u}{\partial\eta} \right) \nonumber \\
&= \frac{d\xi}{dt}\frac{\partial}{\partial\xi}\left( \frac{\partial u}{\partial\xi} + \epsilon\frac{\partial u}{\partial\eta} \right) + \frac{d\eta}{dt}\frac{\partial}{\partial\eta}\left( \frac{\partial u}{\partial\xi} + \epsilon\frac{\partial u}{\partial\eta} \right) \nonumber \\
&= \frac{\partial^2 u}{\partial\xi^2} + 2\epsilon\frac{\partial^2 u}{\partial\xi\partial\eta} + \epsilon^2\frac{\partial^2 u}{\partial\eta^2} \label{ddu2}
\end{align}
Additionally, we have that 
 \begin{equation}
 u(t-\Delta) = u(\xi-\Delta,\eta-\epsilon\Delta)
 \label{hoo1}
 \end{equation}
In discussing the dynamics of \ref{cdddetv}, we will detune the delay $\Delta$ off of its critical value:
 \begin{equation}
\Delta = \Delta_{cr} + \epsilon \Delta_1 + O(\epsilon^2)
\label{hoo2}
 \end{equation}
Substituting Equation \ref{hoo2} into Equation \ref{hoo1} and expanding as a series in $\epsilon$, we get
 \begin{equation}
u(t-\Delta)=\bar u-\epsilon\Delta_1\frac{\partial \bar u}{\partial \xi}-\epsilon\Delta_{cr}\frac{\partial \bar u}{\partial \eta}+O(\epsilon^2)
 \label{hoo3}
 \end{equation}
where $$\bar u \equiv u(\xi-\Delta_{cr},\eta).$$
Now we expand $u$ in a power series in terms $\epsilon$:
 \begin{equation}
u = u_{0} + \epsilon u_{1} + O(\epsilon^2)
 \label{hoo5}
 \end{equation}
 
Substituting Equations \ref{udot2}, \ref{ddu2}, \ref{hoo3} and \ref{hoo5} into Equation \ref{ddot-u}, expanding as a series in $\epsilon$, collecting like terms, and equating like powers of $\epsilon$, we get
 
\begin{align}
 \frac{\partial^2 u_0}{\partial\xi^2} + \frac{\partial u_0}{\partial\xi} + \frac{\lambda}{2\Delta_{cr}}\left( u_0 -\bar u_0 \right) &= 0 \label{e0terms}\\
 \frac{\partial^2 u_1}{\partial\xi^2} + \frac{\partial u_1}{\partial\xi} + \frac{\lambda}{2\Delta_{cr}}\left( u_1 -\bar u_1 \right) &= \left( \alpha  \gamma  \mu  \cos (\gamma  \xi )+\frac{\lambda  \left(\Delta_1-\alpha \Delta_{cr} \sin (\gamma  \xi )\right)}{2 \Delta_{cr}^2} \right)u_0 \nonumber \\
 &\quad + \frac{\lambda  \left(\alpha  \Delta_{cr} \sin (\gamma  \xi )-\Delta_1\right)}{2 \Delta_{cr}^2} \bar u_0 - \mu \frac{\partial u_0}{\partial\eta} -\frac{\lambda}{2}\frac{\partial\bar u_0}{\partial \eta} \nonumber \\
 &\quad + \alpha  \gamma  \cos (\gamma  \xi )\frac{\partial u_0}{\partial \xi} -\frac{\lambda\Delta_1}{2\Delta_{cr}}\frac{\partial\bar u_0}{\partial\xi} - 2\frac{\partial^2 u_0}{\partial\xi\partial\eta} \label{e1terms}
\end{align}

 Equation \ref{e0terms} is linear, constant-coefficient, homogeneous, and does not involve any derivatives with respect to $\eta$. In fact, it is the equation that results from converting the $\epsilon=0$ system (Equations \ref{u-lin0}-\ref{w-lin0}) to a single second-order equation. So we write down the solution:
 \begin{equation}
 u_0 = A(\eta)\cos(\omega_{cr}\xi) + B(\eta)\sin(\omega_{cr}\xi) \label{u0sol}
 \end{equation}
The functions $A(\eta)$ and $B(\eta)$ give the slow flow of the system. We find differential equations on $A(\eta)$ and $B(\eta)$ by substituting Equation \ref{u0sol} into \ref{e1terms} and eliminating the resonant terms. We compute $\bar u_0$ and its partial derivatives using expressions for $\cos(\Delta_{cr}\omega_{cr})$ and $\sin(\Delta_{cr}\omega_{cr})$ given by Equations \ref{crit-2} and \ref{crit-2a}. For example:
\begin{align}
  \bar u_0 &= A(\eta)\cos(\omega_{cr}(\xi-\Delta_{cr})) + B(\eta)\sin(\omega_{cr}(\xi-\Delta_{cr})) \nonumber \\
  &= \left( A(\eta ) \cos (\Delta_{cr} \omega_{cr} )-B(\eta ) \sin (\Delta_{cr} \omega_{cr} ) \right)\cos(\omega_{cr}\xi) \nonumber \\
  &\quad + \left( A(\eta ) \sin (\Delta_{cr} \omega_{cr} )+B(\eta ) \cos (\Delta_{cr} \omega_{cr} ) \right)\sin(\omega_{cr}\xi) \nonumber \\
  &= \left(  B(\eta )\frac{2 \Delta_{cr} \mu }{\lambda } \sqrt{\frac{\lambda -\Delta_{cr} \mu ^2}{\Delta_{cr}}}-A(\eta )\frac{ \left(\lambda -2 \Delta_{cr} \mu ^2\right)}{\lambda } \right)\cos(\omega_{cr}\xi) \nonumber \\
  &\quad + \left(  -A(\eta )\frac{2 \Delta_{cr} \mu }{\lambda } \sqrt{\frac{\lambda -\Delta_{cr} \mu ^2}{\Delta_{cr}}}-B(\eta )\frac{ \left(\lambda -2 \Delta_{cr} \mu ^2\right)}{\lambda } \right)\sin(\omega_{cr}\xi)
\end{align}
  After substituting these expressions into Equation \ref{e1terms} and using angle-sum identities, the remaining trigonometric terms are of the forms
  \begin{equation}
  \cos(\omega_{cr}\xi),\quad \sin(\omega_{cr}\xi),\quad \cos((\omega_{cr}\pm\gamma)\xi),\quad \sin((\omega_{cr}\pm\gamma)\xi)
  \end{equation}
 Notice that $\gamma=2\omega_{cr}$ is the only resonant frequency for the arrival function. For any other value of $\gamma$, the terms involving $\gamma$ at $O(\epsilon)$ are non-resonant, and the first-order two-variable expansion method does not capture any effect from the time-varying arrival function. Therefore, we set $\gamma=2\omega_{cr}$, and Equation \ref{e1terms} becomes
 \begin{align}
   \frac{\partial^2 u_1}{\partial\xi^2} + \frac{\partial u_1}{\partial\xi} + \frac{\lambda}{2\Delta_{cr}}\left( u_1 -\bar u_1 \right) &= \left[c_1 A'(\eta) + c_2 B'(\eta) + c_3 A(\eta) + c_4 B(\eta)\right]\cos(\omega_{cr}\xi) \nonumber \\
   &\quad + \left[c_5 A'(\eta) + c_6 B'(\eta) + c_7 A(\eta) + c_8 B(\eta)\right]\sin(\omega_{cr}\xi) \nonumber \\
   &\quad +\text{ non-resonant terms} \label{u1res}
 \end{align}
 where the coefficients $c_i$ depend on $\lambda,\mu,\alpha,\Delta_{cr}$ and $\Delta_1$.
 
 Elimination of secular terms in Equation \ref{u1res} gives the slow flow equations on $A(\eta)$ and $B(\eta)$:
 \begin{align}
   \frac{dA}{d\eta} &= K_1 A(\eta) + K_2 B(\eta) \\
   \frac{dB}{d\eta} &= K_3 A(\eta) + K_4 B(\eta)
 \end{align}
  where
  \begin{align}
    K_1 &=     \frac{\alpha  \Delta _{cr} \sqrt{\frac{\lambda }{\Delta _{cr}}-\mu ^2} \left(\mu  \Delta _{cr} \left(-4 \mu ^2 \Delta _{cr}+3 \lambda -6 \mu \right)+4 \lambda \right)-2 \Delta _1 \left(\lambda -\mu ^2 \Delta _{cr}\right) \left(\lambda -2 \mu  \left(\mu  \Delta _{cr}+1\right)\right)}{\Delta _{cr} \left(\Delta _{cr} \left(8 \mu ^3 \Delta _{cr}-\lambda ^2-12 \lambda  \mu +12 \mu ^2\right)-16 \lambda \right)} \label{k1val}\\
    K_2 &= \frac{\alpha  \Delta _{cr} \left(\lambda -\mu ^2 \Delta _{cr}\right) \left(-4 \mu ^2 \Delta _{cr}+\lambda -6 \mu \right)+\Delta _1 \sqrt{\frac{\lambda }{\Delta _{cr}}-\mu ^2} \left(\Delta _{cr} \left(-4 \mu ^3 \Delta _{cr}+\lambda ^2+8 \lambda  \mu -4 \mu ^2\right)+8 \lambda \right)}{\Delta _{cr} \left(\Delta _{cr} \left(8 \mu ^3 \Delta _{cr}-\lambda ^2-12 \lambda  \mu +12 \mu ^2\right)-16 \lambda \right)} \\
    K_3 &= \frac{\alpha  \Delta _{cr} \left(\lambda -\mu ^2 \Delta _{cr}\right) \left(-4 \mu ^2 \Delta _{cr}+\lambda -6 \mu \right)+\Delta _1 \sqrt{\frac{\lambda }{\Delta _{cr}}-\mu ^2} \left(\Delta _{cr} \left(4 \mu ^3 \Delta _{cr}-\lambda ^2-8 \lambda  \mu +4 \mu ^2\right)-8 \lambda \right)}{\Delta _{cr} \left(\Delta _{cr} \left(8 \mu ^3 \Delta _{cr}-\lambda ^2-12 \lambda  \mu +12 \mu ^2\right)-16 \lambda \right)} \\
    K_4 &= \frac{\alpha  \Delta _{cr} \sqrt{\frac{\lambda }{\Delta _{cr}}-\mu ^2} \left(\mu  \Delta _{cr} \left(4 \mu ^2 \Delta _{cr}-3 \lambda +6 \mu \right)-4 \lambda \right)-2 \Delta _1 \left(-2 \mu ^2 \Delta _{cr}+\lambda -2 \mu \right) \left(\lambda -\mu ^2 \Delta _{cr}\right)}{\Delta _{cr} \left(\Delta _{cr} \left(8 \mu ^3 \Delta _{cr}-\lambda ^2-12 \lambda  \mu +12 \mu ^2\right)-16 \lambda \right)} \label{k4val}
  \end{align}
  
 The equilibrium point $A(\eta)=B(\eta)=0$ of the slow flow corresponds to a periodic solution for $u_0$, and the stability of the equilibrium corresponds to the stability of that periodic solution. The stability is determined by the eigenvalues of the matrix
  \begin{equation}
  K=
 \left[
 \begin{array}{ccc}
 K_1 & K_2\\
 K_3 & K_4
 \end{array}
 \right]
 \label{mat1}
 \end{equation}
 If both eigenvalues have negative real part, the equilibrium is stable.
 Since the eigenvalues are cumbersome to work with directly, we use the Routh-Hurwitz stability criterion: 
 
 Denote the characteristic polynomial of $K$ by 
 \begin{equation}
 \det(K-r I) = a_0 + a_1 r + a_2 r^2 = 0
 \end{equation}
 Then both eigenvalues have negative real part if and only if all the coefficients satisfy $a_i > 0$. From Equations \ref{k1val}-\ref{k4val}, we have
 \begin{align}
   a_0 &= -\frac{\lambda  \left(\lambda -\Delta _{cr} \mu ^2\right) \left(\Delta_1^2 (\Delta _{cr} (\lambda +4 \mu )+4)-\alpha ^2 \Delta _{cr}^2\right)}{\Delta _{cr}^3 \left(-\Delta _{cr} \lambda ^2-4 \lambda  (3 \Delta _{cr} \mu +4)+4 \Delta _{cr} \mu ^2 (2 \Delta _{cr} \mu +3)\right)} \\
   a_1 &= \frac{4 \Delta_1 \left(\lambda -\Delta _{cr} \mu ^2\right) (\lambda -2 \mu  (\Delta _{cr} \mu +1))}{\Delta _{cr} \left(-\Delta _{cr}\lambda ^2-4 \lambda  (3 \Delta _{cr} \mu +4)+4 \Delta _{cr} \mu ^2 (2 \Delta _{cr} \mu +3)\right)} \\
   a_2 &= 1
 \end{align}
 Recall from Equation \ref{omegadef} that $\omega$ is only positive when $0<\Delta_{cr}<\lambda/\mu^2$. Using this restriction, we find that the coefficients  are all positive when
 
 \begin{equation} \label{tongue1}
 \begin{dcases}
 0 < \lambda \leq 2\mu \\
 \Delta_1>\sqrt{\frac{\alpha ^2 \Delta _{cr}^2}{\Delta _{cr} \lambda +4 \Delta _{cr} \mu +4}}
 \end{dcases}
 \end{equation}
 or alternatively when
  \begin{equation} \label{tongue2}
 \begin{dcases}
  \lambda > 2\mu \\
 0 < \Delta_{cr} < \frac{\lambda -2 \mu }{2 \mu ^2} \\
 \Delta_1 < - \sqrt{\frac{\alpha ^2 \Delta _{cr}^2}{\Delta _{cr} \lambda +4 \Delta _{cr} \mu +4}}
 \end{dcases}
 \end{equation}
 or when
 \begin{equation} \label{tongue3}
 \begin{dcases}
  \lambda > 2\mu \\
 \frac{\lambda -2 \mu }{2 \mu ^2}<\Delta_{cr}<\frac{\lambda }{\mu ^2} \\
 \Delta_1 > \sqrt{\frac{\alpha ^2 \Delta _{cr}^2}{\Delta _{cr} \lambda +4 \Delta _{cr} \mu +4}}
 \end{dcases}
 \end{equation}
 Thus the change of stability occurs at 
 \begin{equation}
 \Delta = \Delta_{cr} \pm \epsilon \sqrt{\frac{\alpha ^2 \Delta _{cr}^2}{\Delta _{cr} \lambda +4 \Delta _{cr} \mu +4}}
 \end{equation}
 where the sign of the $\epsilon$ term depends on $\Delta_{cr}$, $\lambda$ and $\mu$ as in Equations \ref{tongue1}-\ref{tongue3}.
 
\end{proof}
\end{theorem}

\subsection{Numerics for Moving Average Queueing Model}
\begin{figure}
\captionsetup{justification=centering}
		\hspace{-.35in}~\includegraphics[scale = .22]{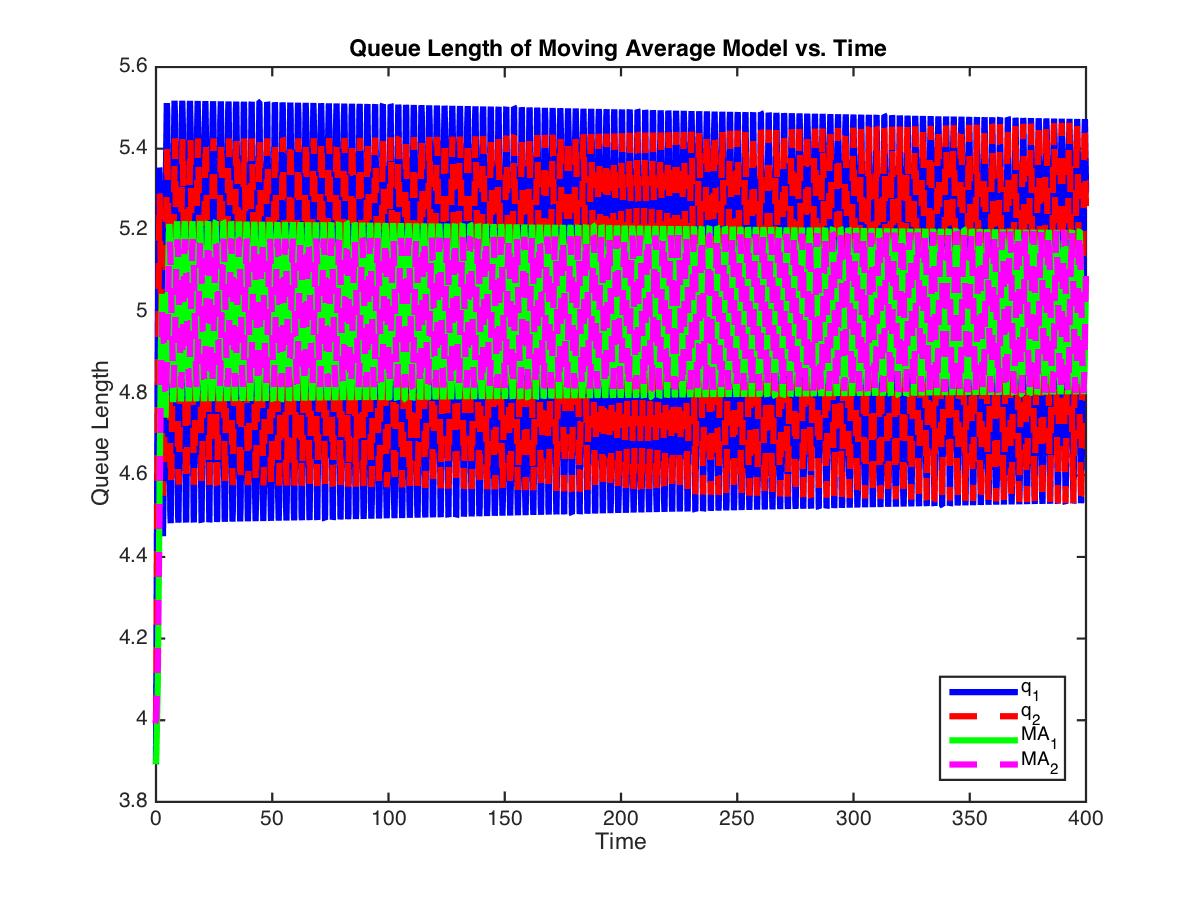}~\hspace{-.4in}~\includegraphics[scale = .22]{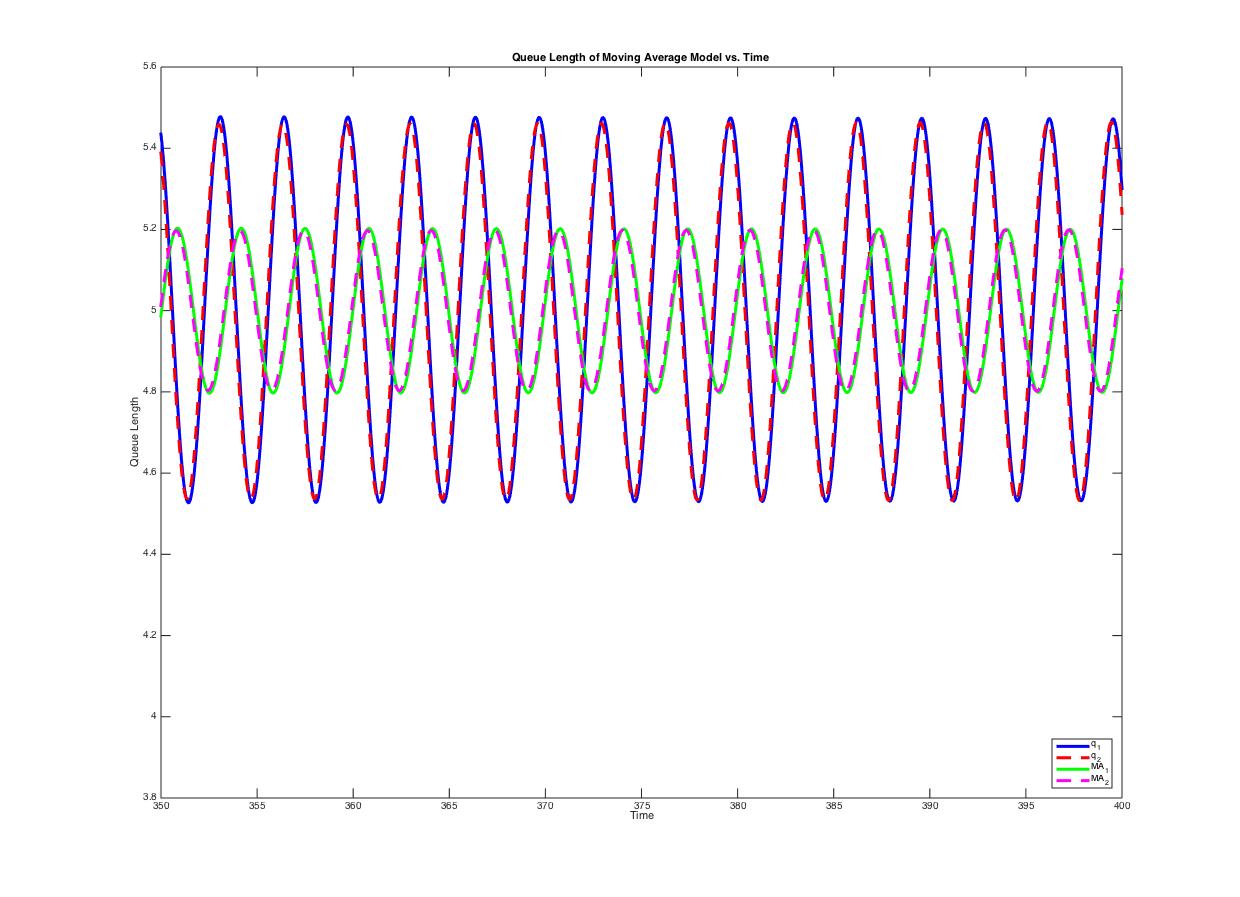}
\caption{ $\Delta_{cr}$ = 2.1448, $\Delta_{mod}$  = 2.2183. \\
$\lambda = 10$, $\mu =1$, $\alpha =1$, $\epsilon = .2$, $\gamma = \sqrt{10/\Delta_{cr} - 1}$, $\Delta$ = 2.18, $\varphi_1([-\Delta,0])=3$,$\varphi_2([-\Delta,0])=4$ } \label{Fig10}
\end{figure}

\begin{figure}
\captionsetup{justification=centering}
		\hspace{-.35in}~\includegraphics[scale = .22]{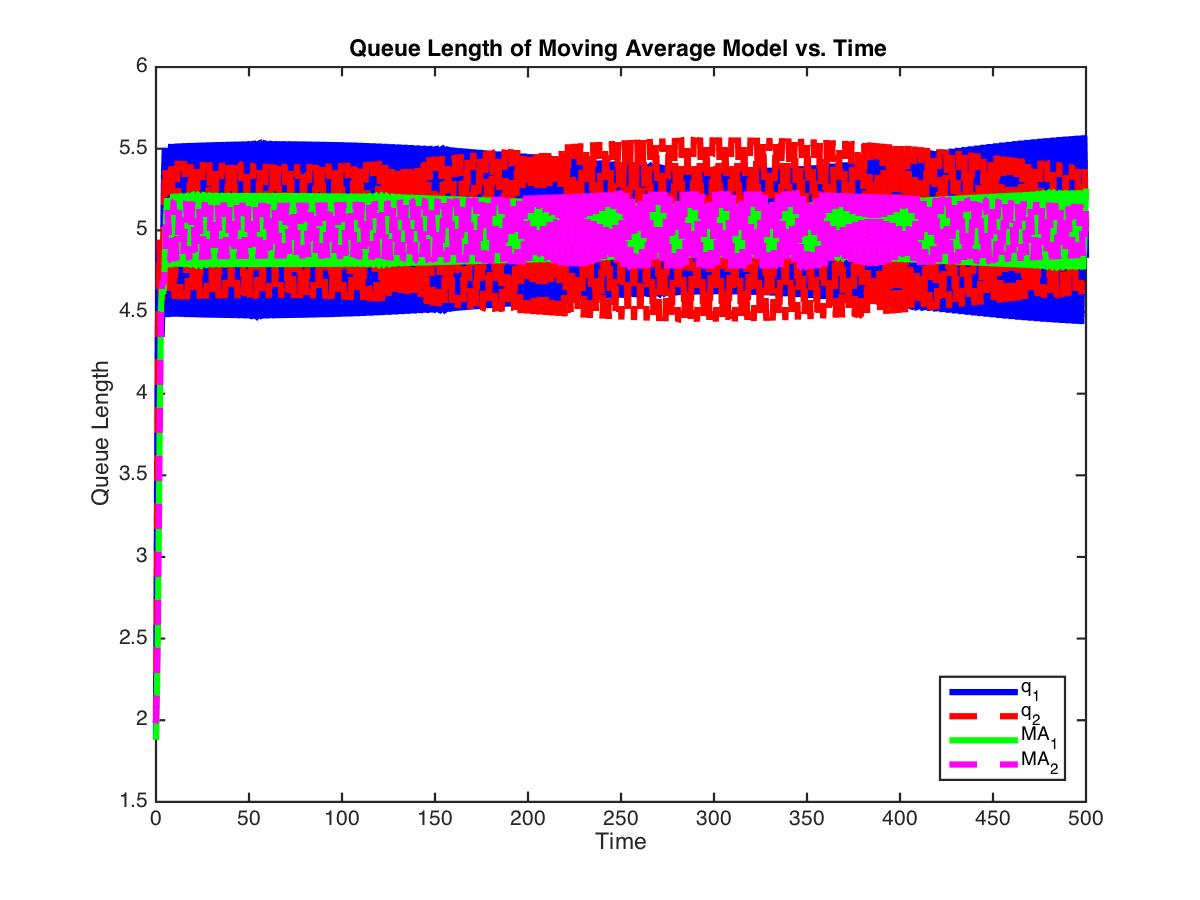}~\hspace{-.4in}~\includegraphics[scale = .22]{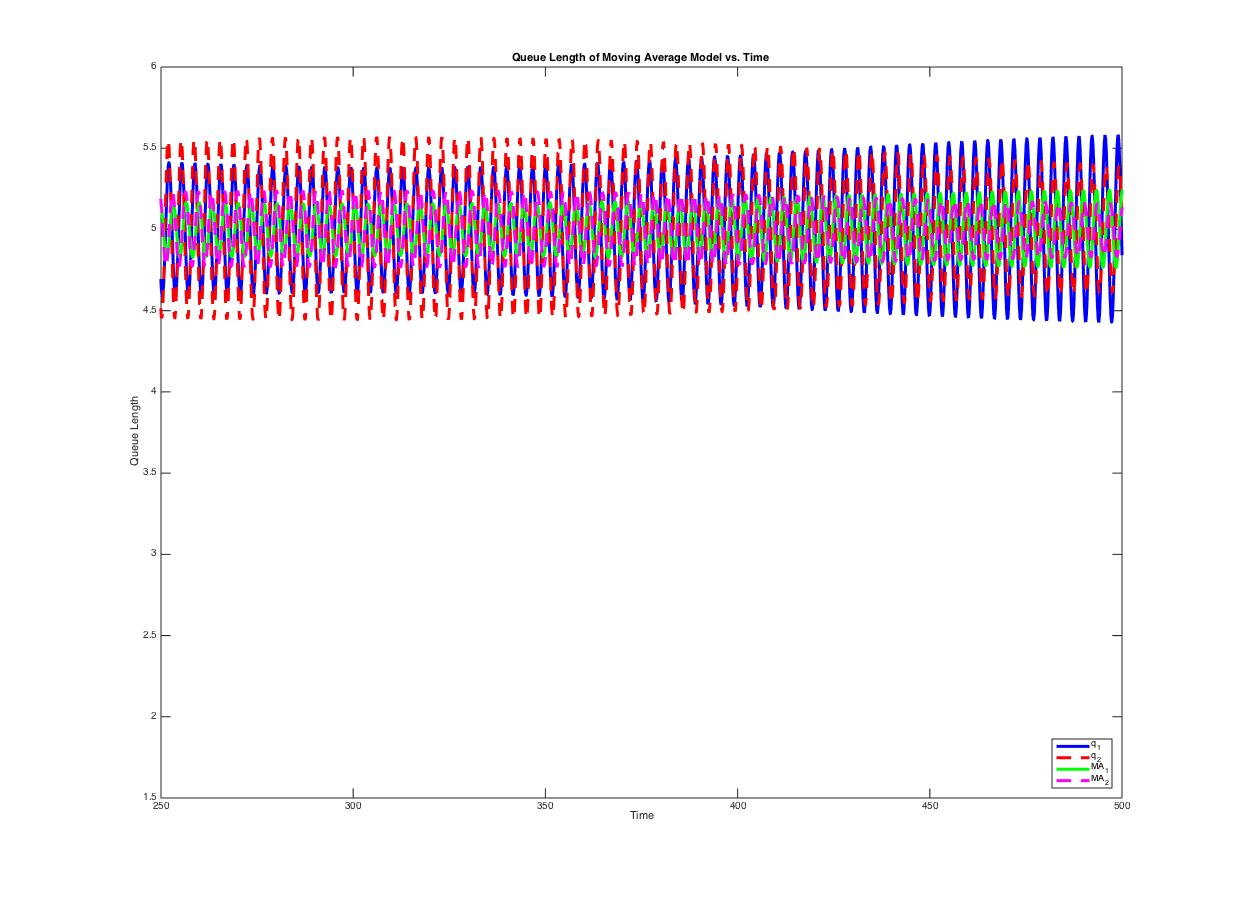}
\caption{ $\Delta_{cr}$ = 2.1448, $\Delta_{mod}$  = 2.2183. \\
$\lambda = 10$, $\mu =1$, $\alpha =1$, $\epsilon = .2$, $\gamma = \sqrt{10/\Delta_{cr} - 1}$, $\Delta$ = 2.25, $\varphi_1([-\Delta,0])=3.9$, $\varphi_2([-\Delta,0])=4$ } \label{Fig11}
\end{figure}

 In this section, we numerically integrate the delay two examples of delay differential equations with moving averages and compare the asymptotic results for determining the Hopf bifurcation that occurs.  On the left of Figure \ref{Fig10} we numerically integrate the two queues and plot the queue length as a function of time.  In this example our lag in information is given by $\Delta =2.18$.  We see that the two equations eventually converge to the same limit as time is increased towards infinity.  This implies that the system is stable and no oscillations or asynchrous dynamics will occur due to instability in this case.  On the right of Figure \ref{Fig10} is a zoomed in version of the figure on the left.   It is clear that the two delay equations are converging towards one another and this system is stable.  However, in Figure \ref{Fig11}. we use the same parameters, but we make the lag in information $\Delta = 2.25$.  This is below the critcal delay in the constant case and above the modified critical delay when the time varying arrival rate is taken into account.  On the right of Figure \ref{Fig11}, we display a zoomed in version of the figure on the left. We see that in this case the two queues oscillate and asynchronous behavior is observed.  Thus, the asymptotic analysis performed works well at predicting the change in stability.    

%\begin{figure}
%\captionsetup{justification=centering}
%		\hspace{-.35in}~\includegraphics[scale = .22]{TV_Paper_MA_1.jpg}~\hspace{-.4in}~\includegraphics[scale = .22]{Closeup_MA_1.jpg}
%\caption{ $\Delta_{cr}$ = 2.1448, $\Delta_{mod}$  = 2.2183. \\
%$\lambda = 10$, $\mu =1$, $\alpha =1$, $\epsilon = .2$, $\gamma = \sqrt{10/\Delta_{cr} - 1}$, $\Delta$ = 2.15, $\varphi_1([-\Delta,0])=1$,$\varphi_2([-\Delta,0])=2$ } \label{Fig9}
%\end{figure}

%\begin{figure}
%\captionsetup{justification=centering}
%		\hspace{-.35in}~\includegraphics[scale = .22]{TV_Paper_MA_3.jpg}~\hspace{-.4in}~\includegraphics[scale = .22]{Closeup_MA_3.jpg}
%\caption{ $\Delta_{cr}$ = 2.1448, $\Delta_{mod}$  = 2.2183. \\
%$\lambda = 10$, $\mu =1$, $\alpha =1$, $\epsilon = .2$, $\gamma = \sqrt{10/\Delta_{cr} - 1}$, $\Delta$ = 2.3, $\varphi_1([-\Delta,0])=1$,$\varphi_2([-\Delta,0])=2$ } \label{Fig12}
%\end{figure}
  
%**************************************************************************
%**************************************************************************

\section{Conclusion and Future Research} \label{sec_conclusion}

 In this paper, we analyze two new two-dimensional fluid models that incorporate customer choice, delayed queue length information, and time varying arrival rates.  The first model considers the customer choice as a multinomial logit model where the queue length information given to the customer is delayed by a constant $\Delta$.  In the second model, we consider customer choice as a multinomial logit model where the queue length information given to the customer is a moving average over an interval of $\Delta$. In the constant arrival case for both models, it is possible to derive an explicit threshold for the critical delay where below the threshold the two queues are balanced and converge to the equilibrium.  However, when the arrival rate is time varying, this problem of finding the threshold is more difficult.  When the time variation is small, we show using asymptotic techniques that the new threshold changes when the arrival rate frequency is twice that of the critical delay frequency. It is important for operators of queues to determine and know these thresholds since using delayed information can have such a large impact on the dynamics of the business.  

Since our analysis is the first of its kind in the queueing literature, there are many extensions that are worthy of future study.  One extension that we would like to explore is the  use of different customer choice functions and incorporating customer preferences in the model.  With regard to customer preferences, this is non-trivial problem because the equilibrium solution is no longer a simple expression, but the solution to a transcendental equation.  This presents new challenges for deriving analytical formulas that determine synchronous or asynchronous dynamics.  A detailed analysis of these extensions will provide a better understanding of what types of information and how the information that operations managers provide to their customers will affect the dynamics of the system.  However, we might be able to use asymptotic techniques for these extensions if we expand around a solution that we know well.  We plan to explore these extensions in subsequent work.

%**************************************************************************
%**************************************************************************

\bibliographystyle{plainnat}
\bibliography{choicetime}
\end{document}